\documentclass{amsart}

\usepackage{amssymb}
\usepackage{eepic}
\usepackage{color}

\allowdisplaybreaks

\numberwithin{equation}{section}

\newtheorem{theorem}{Theorem}[section]
\newtheorem{lemma}[theorem]{Lemma}

\newtheorem{remark}[theorem]{Remark}

\newtheorem{TheoA}{Theorem A}
\newtheorem{TheoB}{Theorem B}
\newtheorem{TheoC}{Corollary C}

\newcommand{\Z}{\mathbb{Z}}
\newcommand{\R}{\mathbb{R}}
\newcommand{\C}{\mathbb{C}}

\newcommand{\summ}{\sum\nolimits}

\def\G{\mathrm{G}}
\def\1{\mathbf{1}}
\def\H{\mathcal{H}}

\def\O{\mathcal{O}}
\def\V{\mathrm{\mathcal{L}(G)}}

\newcommand{\dem}{\noindent {\bf Proof. }}
\newcommand{\demA}{\noindent {\bf End of the proof of Theorem A. }}
\newcommand{\demAA}{\noindent {\bf Proof of ii) $\Rightarrow$ iii) in Theorem A. }}

\newcommand{\demB}{\noindent {\bf Proof of Theorem B. }}
\newcommand{\demC}{\noindent {\bf Proof of Corollary C. }}

\newcommand{\fin}{\hspace*{\fill} $\square$ \vskip0.2cm}

\addtocounter{tocdepth}{-1}

\begin{document}

\title[Twisted Hilbert transforms]{Twisted Hilbert transforms \\ vs Kakeya sets of directions}

\author[Parcet and Rogers]
{Javier Parcet and Keith M. Rogers}




\maketitle



\begin{abstract}
Given a discrete group $\G$ and an orthogonal action $\gamma: \G \to O(n)$ we study $L_p$ convergence of Fourier integrals which are frequency supported on the semidirect product $\R^n \rtimes_\gamma \G$. Given a unit $u \in \R^n$ and $1 < p \neq 2 < \infty$, our main result shows that the twisted (directional) Hilbert transform $H_u \rtimes_\gamma id_\G$ is $L_p$-bounded iff the orbit $\mathcal{O}_\gamma(u)$ is finite. This is in sharp contrast with twisted Riesz transforms $R_u \rtimes_\gamma id_\G$, which are always bounded. Our result characterizes Fourier summability in $L_p$ for this class of groups. We also extend de Leeuw's compactification theorem to this setting and obtain stronger estimates for functions with \lq\lq lacunary" frequency support.
\end{abstract}




\addtolength{\parskip}{+1ex}

\null

\vskip10pt

\null

\section*{Introduction}

Given a $p$-integrable function $f: \R^n \to \C$ with $p > 1$ and a family of bounded functions $m_R$ which converge pointwise to $1$, the classical $L_p$ convergence problem for Fourier integrals consists in determining for which families of $m_R$'s do we find a vanishing limit $$\lim_{R \to \infty} \int_{\R^n} \Big| f(x) - \int_{\R^n} m_R(\xi) \widehat{f}(\xi) \exp(2 \pi i \langle \xi, x \rangle) \, d\xi \Big|^p \, dx \, = \, 0.$$ This problem is typically studied for dilations $m_R(\xi) = m(\xi/R)$ of a fixed function $m: \R^n \to \C$. In this case, the convergence problem reduces to deciding when the Fourier multiplier map associated to $m$ turns out to be $L_p$-bounded. If $m$ is smooth enough, the H\"ormander-Mihlin criterium suffices to obtain $L_p$ convergence. The problem is more interesting when the multiplier is less smooth, like $m = \chi_\Omega$ for some (convex) open set $\Omega$ containing the origin. The $L_p$-boundedness of directional Hilbert transforms shows immediately that $\Omega$ can always be taken to be a convex polyhedron. On the other hand, a consequence of Fefferman's multiplier theorem for the ball \cite{Fef} is that the boundary of $\Omega$ must be flat. These considerations also apply to periodic functions $f: \mathbb{T}^n \to \C$ if we replace Fourier integrals by Fourier series, with  frequency group given by the integer lattice $\Z^n$. Moreover, according to de Leeuw's theorem \cite{dL}, the same discussion is valid for the Bohr compactification of $\R^n$, whose frequency group is $\R^n$ with the discrete topology $\R^n_{\mathrm{disc}}$. What happens for other frequency groups? We will give necessary and sufficient conditions on convex polyhedra so that there is Fourier summability in $L_p$ for semidirect products of the classical groups above with arbitrary discrete groups acting on them. Other classes of groups will be also considered.

The group von Neumann algebra associated to a locally compact group is a noncommutative analog of the algebra of essentially bounded functions on the dual group of a given abelian group. As basic models of quantum groups, they play a prominent role in noncommutative geometry and operator algebra \cite{Con,KR}. Harmonic analysis on these algebras puts the reference group on the frequency side, in contrast with the vast literature generalizing classical groups on the spatial side, like in the case of nilpotent Lie groups \cite{M,M2,M3,MPR,MRS}. Our dual approach is inspired by the ground-breaking results of Cowling/Haagerup \cite{CH,Haa} on the approximation property and Fourier multipliers on group von Neumann algebras. This paper is part of an effort to extend Fourier analysis to discrete group von Neumann algebras \cite{Boz,H,JM,JMP1,JMP3}. We will study directional Hilbert transforms (semispace Fourier multipliers) and the $L_p$ convergence of Fourier series in this setting. This is a natural continuation of \cite{JM,JMP1}, where smooth multipliers were investigated. We will work with semidirect products $\R^n \rtimes \G$, $\R^n_{\mathrm{disc}} \rtimes \G$ and $\Z^n \rtimes \G$ with an arbitrary discrete group $\G$, which are basic models of nonabelian groups which contain a Euclidean subgroup embedded in a nontrivial way. This class of groups is sufficiently rich to exhibit barriers to Fourier summability not encountered in the classical case. It includes well-known examples (via subgroups) like the discrete Heisenberg group or the infinite dihedral group, which have been studied so far on the spatial side but not on the frequency side.
Our methods also establish sufficient conditions for $L_p$ convergence on more general nonabelian discrete groups.

Let $\G$ be a discrete group with left regular representation $\lambda_\G: \G \to \mathcal{B}(\ell_2(\G))$ given by $\lambda_\G(g) \delta_h = \delta_{gh}$, where the $\delta_g$'s form the unit vector basis of $\ell_2(\G)$. Write $\mathcal{L}(\G)$ for its group von Neumann algebra, the weak operator closure of the linear span of $\lambda_\G(\G)$. Consider the standard trace $\tau_\G(\lambda_\G(g)) = \delta_{g=e}$ where $e$ denotes the identity of $\G$. Any $f \in \mathcal{L}(\G)$ has a Fourier series $\sum_g \widehat{f}(g) \lambda_\G(g)$ with $\tau_\G(f) = \widehat{f}(e)$ and we set $$L_p(\widehat{\mathbf{G}}) = L_p(\mathcal{L}(\mathrm{G}), \tau_\G) \equiv \mbox{Closure of $\mathcal{L}(\G)$ wrt $\|f\|_{L_p(\widehat{\mathbf{G}})} = \big( \tau_\G [ |f|^p ] \big)^\frac1p$},$$ the natural $L_p$ space over the noncommutative measure space $(\mathcal{L}(\mathrm{G}), \tau_\G)$. We invite the reader to check that $L_p(\mathcal{L}(\Z^n), \tau_{\Z^n}) = L_p(\mathbb{T}^n)$, after identifying $\lambda_{\Z^n}(k)$ with $e^{2\pi i \langle k, \cdot \rangle}$. In general, the (unbounded) operator $|f|^p$ is obtained from functional calculus on the Hilbert space $\ell_2(\G)$. Assume now that $\G$ acts on $\R^n$ by orthogonal transformations and let $$\gamma: \G \to O(n)$$ stand for the corresponding action. Let us write $\R_{\mathrm{disc}}^n$ for the $n$-dimensional Euclidean space equipped with the discrete topology. As a discrete abelian group we find $\mathcal{L}(\R_{\mathrm{disc}}^n) \simeq L_\infty(\R_{\mathrm{bohr}}^n)$, the algebra of essentially bounded functions on the Bohr compactification. The semidirect product $\Gamma_{\mathrm{disc}} = \R_{\mathrm{disc}}^n \rtimes_\gamma \G$ is still discrete and elements of $\mathcal{L}(\Gamma_{\mathrm{disc}})$ are formally given by $$f \sim \sum_{\xi \in \R_{\mathrm{disc}}^n} \sum_{g \in \G} \widehat{f} (\xi,g) \lambda_{\Gamma_{\mathrm{disc}}}(\xi \rtimes_\gamma g) \sim \sum_{g \in \G} f_g \rtimes_\gamma \lambda_\G(g)$$ with $f_g \sim \sum_\xi \widehat{f}(\xi,g) \mbox{b-exp}_\xi$ and $\mbox{b-exp}_\xi$ the $\xi$-th character on $\R^n_{\mathrm{bohr}}$. Its restriction to $\R^n$ is the standard character $\exp_\xi(x) = \exp(2 \pi i \langle \xi, x \rangle)$. The formal equivalence follows from the isometric isomorphism $\mathcal{L}(\Gamma_{\mathrm{disc}}) \simeq \mathcal{L}(\R_{\mathrm{disc}}^n) \rtimes_\gamma \G$ with the cross product algebra, whose main operations are recalled below:
\begin{itemize}
\item $(f \rtimes_\gamma \lambda_\G(g))^* = \gamma_{g^{-1}}(f^*) \rtimes_\gamma \lambda_\G(g^{-1})$,

\vskip6pt

\item $(f \rtimes_\gamma \lambda_\G(g)) (f' \rtimes_\gamma \lambda_\G(g')) = f \gamma_g(f') \rtimes_\gamma \lambda_\G(gg')$,

\item $\tau_{\mathcal{L}(\R_{\mathrm{disc}}^n)} \rtimes_\gamma \tau_\G(f \rtimes_\gamma \lambda_\G(g)) = \displaystyle \delta_{g=e} \int_{\R_{\mathrm{bohr}}^n} f(x) d\mu(x)$,
\end{itemize}
with $\gamma_gf(x) = f(\gamma_{g^{-1}}x)$ and $\mu$ the normalized Haar measure on $\R_{\mathrm{bohr}}^n$. Restoring the usual topology on $\R^n$, the algebra $\mathcal{L}(\Gamma)$ for $\Gamma = \R^n \rtimes_\gamma \G$ can still be represented as $\mathcal{L}(\R^n) \rtimes_\gamma \G$, where $\mathcal{L}(\R^n) \simeq L_\infty(\R^n)$ with the Lebesgue measure and formally we have $$f \sim \int_{\R^n} \sum_{g \in \G} \widehat{f}(\xi,g) \lambda_{\Gamma}(\xi \rtimes_\gamma g) \, d\xi.$$ Given a convex polyhedron $\mathrm{K}$ containing the origin of $\R^n$, the question is whether we still have $L_p$ convergence of Fourier series/integrals along dilations of $\mathrm{K}$. In other words, under which conditions do we have
\begin{itemize}
\item[Q1)] $\displaystyle \lim_{R \to \infty} \Big\| f - \int_{R \mathrm{K}} \, \sum_{g \in \G} \widehat{f}(\xi,g) \lambda_{\Gamma}(\xi \rtimes_\gamma g) \, d\xi \, \Big\|_{L_p(\widehat{\mathbf{\Gamma}})} \hskip0.5pt = \, 0$,

\item[Q2)] $\displaystyle \lim_{R \to \infty} \Big\| f - \sum_{\xi \in R\mathrm{K}} \sum_{g \in \G} \widehat{f} (\xi,g) \lambda_{\Gamma_{\mathrm{disc}}}(\xi \rtimes_\gamma g) \Big\|_{L_p(\widehat{\mathbf{\Gamma}_{\mathrm{disc}}})} = \, 0$,

\vskip5pt

\item[Q3)] $L_p$ convergence of Fourier series for more general groups.
\end{itemize}
According to de Leeuw's theorem \cite{dL}, the $L_p(\R^n_{\mathrm{bohr}})$-boundedness of $$H_u: \summ_\xi \widehat{f}(\xi) \, \mbox{b-exp}_\xi \mapsto - i \summ_\xi \mathrm{sgn} \langle \xi,u \rangle \widehat{f}(\xi) \, \mbox{b-exp}_\xi \qquad (u \in \mathbb{S}^{n-1})$$ is equivalent to the $L_p$-boundedness of $\widehat{H_uf}(\xi) = - i \mathrm{sgn} \langle \xi,u \rangle \widehat{f}(\xi)$, the $u$-directional Hilbert transform in $\R^n$ with the usual topology. The map $\exp_k \mapsto - i \hskip1pt \mathrm{sgn} \langle k,u \rangle \exp_k$ is also $L_p$ bounded for $\G = \Z^n$ by transference arguments. In the context of Lie groups, we may consider similar operators by means of the exponential map with the vector $u$ living in the corresponding Lie algebra. Motivated by our questions above, we are interested in directional Hilbert transforms for the class of groups $\Gamma$ and $\Gamma_{\mathrm{disc}}$. We characterize the $L_p$-boundedness of $H_u \rtimes_\gamma id_\G$ for $1 < p \neq 2 < \infty$ on both algebras.


\begin{TheoA}
Consider the operator densely defined by $$H_u \rtimes_\gamma id_\G: \summ_g^{\null} f_g \rtimes_\gamma \lambda_\G(g) \, \mapsto \, \summ_g H_u(f_g) \rtimes_\gamma \lambda_\G(g).$$ If $1 < p \neq 2 < \infty$ and $u \in \mathbb{S}^{n-1}$, the following properties are equivalent
\begin{itemize}
\item[i)] $H_u \rtimes_\gamma id_\G$ is bounded on $L_p(\widehat{\mathbf{\Gamma}})$,

\vskip3pt

\item[ii)] $H_u \rtimes_\gamma id_\G$ is bounded on $L_p(\widehat{\mathbf{\Gamma}_{\mathrm{disc}}})$,

\vskip5pt

\item[iii)] The orbit $\mathcal{O}_\gamma(u) = \{ \gamma_g(u) \, | \, g \in \G \}$ is finite,

\vskip5pt

\item[iv)] The following matrix inequality holds $$\int_{\R^n} \Big\| \Big( H_{\gamma_{g^{-1}}(u)}(f_{g,h})(x) \Big) \Big\|_{S_p(\G)}^p \, dx \, \le c_p \, \int_{\R^n} \Big\| \Big( f_{g,h}(x) \Big) \Big\|_{S_p(\G)}^p \, dx.$$
\end{itemize}
We will also prove $L_1 \to L_{1,\infty}$ and $L_\infty \to \mathrm{BMO}$ type estimates for finite orbits.
\end{TheoA}

Our arguments combine Kakeya type constructions, ergodic methods, geometric group theory and noncommutative Littlewood-Paley decompositions. The main application of Theorem A is that we find complete answers to Questions 1,2 and partial answers to Question 3 in Corollary C below. Although we refer (for obvious reasons) to these operators as \emph{$\gamma$-twisted Hilbert transforms}, we have not found a close relation to other twisted convolution operators \cite{GL,M,St2}. Of course, twisted Hilbert transforms are always $L_2$ bounded since they can be realized as Fourier multipliers on $\mathcal{L}(\Gamma)$/$\mathcal{L}(\Gamma_{\mathrm{disc}})$ with an $\ell_\infty$ multiplier. As an extension of $H_u$ to a larger space, is it conceivable that $H_u \rtimes_\gamma id_\G$ should remain bounded in $L_p$ for $1 < p < \infty$. In fact, as recently proven in \cite{JMP1} this is exactly what happens for Riesz transforms $$R_u \rtimes_\gamma id_\G: \summ_g f_g \rtimes \lambda_\G(g) \mapsto \summ_g R_u(f_g) \rtimes \lambda_\G(g).$$ In contrast, Theorem A establishes a surprisingly rigid characterization in terms of the $\gamma$-orbit of $u$. This rigidity has led us to analyze the behavior of $(H_u \rtimes_\gamma id_\G)f$ when the Fourier spectrum of $f$ is supported on a subset $\Lambda$ of $\G$ leading to infinite but lacunary $\gamma$-suborbits of $u$, see Theorem B below. Our notion of lacunarity partly relies on suitable length functions/cocycles on $\G$, emphasizing after \cite{JMP1} the role of cohomology theory in our approach.

We will prove Theorem A by showing i) $\Rightarrow$ ii) $\Rightarrow$ iii) $\Rightarrow$ i). The additional equivalence with iv) and the endpoint estimates will be proven later. The core of the proof is ii) $\Rightarrow$ iii). Roughly, when $\mathcal{O}_\gamma(u)$ is not finite, our strategy is to construct a Littlewood-Paley type decomposition determined by a sequence $g_1, g_2, \ldots$ in a group amplification of $\G$ so that $\gamma_{g_1}(u), \gamma_{g_2}(u), \ldots$ admits Kakeya sets of directions in some sense. The idea is then to show that $L_p$-boundedness of $H_u \rtimes_\gamma id_\G$ implies a \lq twisted Meyer inequality' in the compactified space $L_p(\R^n_\mathrm{bohr})$, $$\Big\| \Big( \sum_{j=1}^\infty |H_{\gamma_{g_j}(u)}(f_{g_j})|^2 \Big)^\frac12 \Big\|_p \, \le \, c_p \max \left\{ \Big\| \Big( \sum_{j=1}^\infty |f_{g_j}|^2 \Big)^\frac12 \Big\|_p, \Big\| \Big( \sum_{j=1}^\infty |\gamma_{g_j^{-1}}f_{g_j}|^2 \Big)^\frac12 \Big\|_p \right\}$$ whenever $p > 2$; a dual inequality arises for $p < 2$. The maximum on the right comes from the row/column nature of noncommutative square functions, \hskip-1pt see \cite{JMP1,MP,PisAJM,PX} for more on Littlewood-Paley estimates in this context. See also Lemma \ref{SqMax} for more on Meyer type inequalities and Remark \ref{RMeyer} for twisted generalizations. A less tractable square function estimate in $L_p(\R^n)$ is obtained by decompactification. The goal is then to pick directions $\gamma_{g_j}(u)$ distributed in a way so that such an inequality cannot occur for a suitable choice of the functions $f_{g_j}$. This is reminiscent of Fefferman's estimates in the solution of the ball multiplier problem \cite{Fef}.
The additional term on the right-hand side is estimated by ergodic arguments. We will finally disprove such an inequality for sets of directions that have a $2$-dimensional projection which admits Besicovitch's sprout triangle construction. Let us say that such a set \emph{admits Kakeya shadows}. Using further group theoretical tools, we will study the geometry of $\mathcal{O}_\gamma(u)$ and show that infinite orbits always admit Kakeya shadows. Our choice of functions $f_{g_j}$ adapts Fefferman's construction to higher dimensions in a way which is compatible with the second (new) term in the maximum considered above. On the other hand, the boundedness of $H_u \rtimes_\gamma id_\G$ in $L_p$ for finite orbits can be proved by standard methods.

The equivalence i) $\Leftrightarrow$ ii) may be regarded as a twisted form of K. de Leeuw's compactification theorem \cite{dL} for the multiplier $H_u \rtimes_\gamma id_\G$. Some other results along this line will appear in \cite{JMP2,PPR} as a byproduct of noncommutative transference methods. According to the sketched argument ---a decompactification process and a Kakeya type set construction--- the implication ii) $\Rightarrow$ i) in the twisted form of de Leeuw's theorem is now reduced to the straightforward one iii) $\Rightarrow$ i). Our interest in ii) $\Leftrightarrow$ iv) was motivated by the Neuwirth/Ricard transference argument \cite{NR} which provides a slightly weaker result for amenable $\G$. Instead, we notice the equivalence iii) $\Leftrightarrow$ iv) for general $\G$, which follows easily from our previous approach. Although not closely related, it is somehow amusing to compare this  with the Bateman/Thiele results \cite{B2,BT} on Hilbert transforms along one-variable vector fields, see also the work of Lacey/Li \cite{LL1,LL2}. They consider operators on $\R^2$ of the form $Tf(x,y) = H_{u(x)}f(x,y)$ so that the directions change only with the first variable. In our case, the operator lives in the algebra of matrix-valued functions and the directions change only with the rows $\langle e_g, Tf(x) e_h \rangle = H_{\gamma_g(u)}f_{g,h}(x)$. In the case of finite orbits, we also provide the associated weak type $L_1$ inequality for amenable groups and certain $L_\infty \to \mathrm{BMO}$ estimate for general discrete groups. The weak type inequality arises from \cite{NR,Pa1}. The latter requires a suitable choice of \emph{$u$-directional} BMO, which might be new even in the Euclidean case.

Our remaining results require some terminology. A set $\Omega = \{ \omega_j \, | \, j \ge 1\}$ in the unit sphere $\mathbb{S}^{n-1}$ will be called \emph{radially lacunary} if there exists a limit point $\omega$ in the sphere such that $$\sup_{j \ge 1} \frac{|\omega_{j+1} - \omega|}{|\omega_j - \omega|} \, < \, 1.$$ This \lq radial' lacunarity is in some sense one-dimensional. In a recent paper \cite{PR} on directional maximal operators, we introduced a higher dimensional notion. Given $\Omega$ as above, set $d = \dim \, [\mathrm{span}(\Omega)]$ and $\Sigma(d) = \{(j,k) \, | \ 1 \le j < k \le d\}$. Divide $\Omega$ into lacunary segments with respect to an orthonormal basis $e_1, e_2, \ldots, e_d$ $$ \Omega_{\sigma\!,i}=\Big\{\, \omega \in \Omega \ \big| \ 0 < \theta_{\sigma\!, i+1} < \Big|\frac{\langle \omega, e_k \rangle}{\langle \omega, e_j \rangle} \Big| \le \theta_{\sigma\!, i}\, \Big\} \quad \mbox{for} \quad \sigma = (j,k) \in \Sigma(d)$$ with
$\sup_{i \in \Z} \theta_{\sigma\!, i+1}/\theta_{\sigma\!,i} < 1$. Letting $\Omega_{\sigma\!,\infty} = \Omega\cap (e_j^\perp \cup e_k^\perp)$ and $\Z^* = \mathbb{Z}\cup\{\infty\}$ we obtain a partition $\{\Omega_{\sigma\!,i}\}_{i\in\mathbb{Z}^*}$ of $\Omega$ for all $\sigma \in \Sigma(d)$. A \emph{dissection} will be such a choice of ${d}\choose{2}$ partitions. The set of directions $\Omega$ is called $\mathrm{HD}$-lacunary of order $0$ if it consists of a single direction. Recursively, it is $\mathrm{HD}$-\emph{lacunary of order} $L$  if there is a dissection for which the $\Omega_{\sigma\!,i}$'s are HD-lacunary of order $\le L-1$ for all $i\in \mathbb{Z}^*$ and
$\sigma \in \Sigma(d)$ with uniformly bounded lacunary constants.

Given $\Lambda \subset \G$, we say that $(\Lambda,\gamma,u)$ is a \emph{lacunary triple} when the $\gamma$-suborbit $\mathcal{O}_\gamma(\Lambda^{-1},u) = \{\gamma_g(u) \, | \, g^{-1} \in \Lambda\}$ is radially lacunary and HD-lacunary of finite order. Note that neither notion of lacunarity is stronger than the other. We also introduce the space
$$L_{\Lambda,p}(\widehat{\mathbf{\Gamma}_{\mathrm{disc}}}) \, = \, \Big\{ f \sim \sum_{g \in \Lambda} f_g \rtimes_\gamma \lambda_\G(g) \in L_p(\widehat{\mathbf{\Gamma}_\mathrm{disc}}) \Big\}.$$

\begin{TheoB} We have $$H_u \rtimes_\gamma id_\G: L_{\Lambda,p}(\widehat{\mathbf{\Gamma}_\mathrm{disc}}) \to L_{\Lambda,p}(\widehat{\mathbf{\Gamma}_\mathrm{disc}})$$ for any $\Lambda \subset \G$ for which $(\Lambda,\gamma,u)$ is a lacunary triple and for any $1 < p < \infty$. 
\end{TheoB}

In other words, this result gives a sufficient condition on $\Lambda$ for $L_p$-boundedness of $\gamma$-twisted Hilbert transforms acting on functions $f = \sum_{g \in \Lambda} f_g \rtimes_\gamma \lambda_\G(g)$ whose Fourier spectrum lies in $\Lambda$. In particular, it provides infinite dimensional subspaces on which $H_u \rtimes_\gamma id_\G$ is bounded when $\mathcal{O}_\gamma(u)$ is not finite. It is quite simple to construct specific examples. According to Lemma \ref{SqMax} and recent estimates for the directional maximal function \cite{PR}, $\mathrm{HD}$-lacunarity is used to avoid the presence of Kakeya shadows and radial lacunarity is used to apply Littlewood-Paley type estimates. The analog of Theorem B with $\Gamma$ in place of $\Gamma_{\mathrm{disc}}$ only requires the adaptation of the Littlewood-Paley estimates in \cite{JMP1} to the (non-discrete) group~$\Gamma$. More general notions of lacunarity provide generalizations of Theorem B, see Remark \ref{GralLac}. A more in depth analysis would be related to some classical problems in harmonic analysis.

We have shown how Lie groups or semidirect products $\R^n \rtimes_\gamma \G$ admit enough geometric structure to define directional Hilbert transforms on them. Given a general discrete group, there is no standard \lq space of directions\rq${}$ to define Hilbert transforms on its group algebra, as we could do with $\R^n$ or the corresponding Lie algebra. According to \cite{JMP1}, the key point is to use a broader interpretation of tangent spaces in terms of length functions and cocycles. These tools provide natural forms of directional Hilbert transforms. Moreover, the problem of $L_p$ convergence for Fourier series can be reformulated for general group von Neumann algebras in terms of cocycles. Given a conditionally negative length function $\psi: \G \to \R_+$ with associated cocycle $(\H, b, \gamma)$ ---precise definitions can be found Section \ref{Sect1}--- consider any open convex bounded polytope $\mathrm{K}$ in the Hilbert space $\H$ containing the origin. Typically we may think of $\mathrm{K}$ as a cube centered at the origin. The problem is then to determine conditions on $\mathrm{K}$, for which truncation along dilations of $\mathrm{K}$ yields $L_p$ convergence of the partial sums; $$\lim_{R \to \infty} \big\| f - T_{R,\psi} f \big\|_{L_p(\widehat{\mathbf{G}})} = 0 \qquad \mbox{where} \qquad T_{R,\psi} f = \sum_{g: \hskip3pt b(g) \in R \mathrm{K}} \widehat{f}(g) \lambda_\G(g).$$ Of course, inner cocycles are less interesting in this regard since partial sums are not finite truncations and the norm limit stabilizes in finite time. In the next result we provide sufficient conditions for $L_p$ convergence on any pair $(\G,\psi)$. Moreover, we give complete answers to Questions 1, 2 and partial answers to Question 3. Given $u \in \H$, define the $(\psi,u)$-Hilbert transform as $$H_{\psi,u}: \lambda_\G(g) \mapsto -i \langle b(g), u \rangle_{\H} \lambda_\G(g).$$ Again, the space $L_{\Lambda,p}(\widehat{\mathbf{G}})$ is the closure in $L_p(\widehat{\mathbf{G}})$ of elements $f \sim \sum_{g \in \Lambda} \widehat{f}(g) \lambda_\G(g)$.

\begin{TheoC}
Let $\G$ be a discrete group equipped with a length function $\psi: \G \to \R_+$ with associated $($finite-dimensional\hskip1pt $)$ cocycle $(\H, b, \gamma)$. If $1 <p < \infty$, then the following results hold\hskip1pt $:$
\begin{itemize}
\item[a)] $H_{\psi,u}: L_p(\widehat{\mathbf{G}}) \to L_p(\widehat{\mathbf{G}})$ if $\mathcal{O}_{\gamma}(u)$ is a finite orbit.

\vskip3pt

\item[b)] $H_{\psi,u}:  L_{\Lambda,p}(\widehat{\mathbf{G}}) \to L_{\Lambda,p}(\widehat{\mathbf{G}})$ if $(\Lambda, \gamma,u)$ is a lacunary triple.

\vskip3pt

\item[c)] $T_{R,\psi} f \to f$ in $L_p(\widehat{\mathbf{G}})$ whenever $|\mathcal{O}_{\gamma}(u_j)| < \infty$ for the normal directions $u_1, u_2, \ldots, u_m$ to all the faces of $\mathrm{K}$. Moreover, this condition is necessary for $L_p$ convergence in $\Gamma = \R^n \rtimes_\gamma \G$ or $\Gamma_{\mathrm{disc}} = \R^n_{\mathrm{disc}} \rtimes_\gamma \G$.
\end{itemize}
\end{TheoC}

Note that, since $\mathrm{K}$ is a convex polyhedron, the directions $u_1, u_2, \ldots, u_m$ must span $\R^n$. In particular, if the $\gamma$-orbit of the $u_j$'s is finite for all $1 \le j \le m$, we must have $\mathrm{N} = \sup_{\xi \in \R^n} |\mathcal{O}_\gamma(\xi)| < \infty$. This might appear to be quite a rigid condition in the line of Theorem A, but we recall that the classical theory has only studied $\Z^n$ and $\R^n$ via the trivial cocycle, which is given by the identity inclusion map with the trivial action, so that orbits in that case are always of cardinality $1$. A simple example for which we may find $1 < \mathrm{N} < \infty$ is given by the infinite dihedral group $\Z_2 * \Z_2 \simeq \Z \rtimes \Z_2$. On the other hand, this strong rigidity is surprisingly necessary for both $\Gamma$ and $\Gamma_{\mathrm{disc}}$. Additionally, for crossed products with the integer lattice $\Z^n \rtimes_\gamma \G$ all possible actions fixing $\Z^n$ must satisfy $\mathrm{N} < \infty$. The proof of Corollary C follows from our previous results on $\Gamma_{\mathrm{disc}}$ together with an intertwining identity from \cite{JMP1}. The condition in a) no longer characterizes $L_p$ boundedness of $H_{\psi,u}$ in general. Indeed, take $\G = \Z$ with $(\H,b,\gamma)$ given by $$\H = \C \simeq \R^2, \quad b(k) = \exp_\alpha(k)-1, \quad \gamma_k(z) = \exp_\alpha(k)z$$ for some $\alpha \in \R_+$. When $\alpha \in \R \setminus \mathbb{Q}$, the orbit of any $u \in \mathbb{S}^1$ is not finite, but the map $H_{\psi,u}: \exp_k \mapsto - i \, \mathrm{sgn} \langle \exp_\alpha(k) - 1, u \rangle \exp_k$ is $L_p$-bounded. This easily follows from de Leeuw's periodization and restriction theorems \cite{dL}.

Corollary C should be compared with the $L_p$ convergence problem of Fourier series/integrals for vector-valued functions, where directional Hilbert transforms are always $L_p$-bounded as long as we take values in a UMD Banach space. If we replace cross products by tensor products in our setting, we obtain $$L_p(\mathcal{L}(\R^n) \otimes \mathcal{L}(\G)) = L_p(\mathcal{L}(\R^n \times \G)) = L_p(\R^n; L_p(\mathcal{L}(\G)))$$ and $L_p$ convergence for the group $\R^n \times \G$ follows from the vector-valued theory since noncommutative $L_p$ spaces are UMD for any $1 < p < \infty$. Our conditions in Corollary~C for $\R^n \rtimes \G$ are (necessarily) much more demanding.

\section{A Littlewood-Paley type theorem} \label{Sect1}

The first ingredient for the proof of Theorem A will be a Littlewood-Paley type inequality for group von Neumann algebras. More concretely, let $\G$ be a discrete group equipped with a given length function $\psi: \G \to \R_+$ and consider a lacunary partition of $\R_+ = \bigcup_m \mathcal{I}_m$. Then, any $f \in L_p(\widehat{\mathbf{G}})$ can essentially be written as $f \sim \sum_m f_m$, where the Fourier spectrum of $f_m$ lies in $\Lambda_m = \{ g \in \G \, \big| \ \psi(g) \in \mathcal{I}_m \}$. One of the main results from \cite{JMP1} provides a norm equivalence between $f$ and certain noncommutative square function associated to the $f_m$'s. In this section, we recall this Littlewood-Paley estimate and consider some natural length functions in $\G$ related to $\gamma$. We also refer to Pisier's papers \cite{PisAJM,P2} for more on lacunary type sets in discrete groups and noncommutative Littlewood-Paley inequalities.

\subsection{Length functions and cocycles}

An affine representation of $\G$ is a group homomorphism $\G \to \H \rtimes O(\H)$ into the affine group associated to a real Hilbert space $\H$. Affine representations are determined by a representation $\gamma: \G \to O(\H)$ together with a mapping $b: \G \to \H$ satisfying $b(gh) = \gamma_g(b(h)) + b(g)$. The triple $(\H, \gamma, b)$ is usually referred to as a cocycle of $\G$. It is clear that affine representations and cocycles of $\G$ are in one-to-one correspondence. In this paper, we say that $\psi: \G \to \R_+$ is a \emph{length function} if it vanishes at the identity $e$, $\psi(g) = \psi(g^{-1})$ and $\summ_g \beta_g = 0 \Rightarrow \summ_{g,h} \overline{\beta}_g \beta_h \psi(g^{-1}h) \le 0$. Those functions satisfying the last condition are called conditionally negative. It is straightforward to show that length functions take values in $\R_+$. Length functions are also in one-to-one correspondence with affine representations and cocycles. Namely, any cocycle $(\H, \gamma,b)$ gives rise to the length function $\psi_b(g) = \langle b(g), b(g) \rangle_\H$, as it can be easily checked. Reciprocally, any length function $\psi$ gives rise to a cocycle $(\H, \gamma, b)$. This is a standard application of Schoenberg's theorem \cite{Sc}, which claims that $\psi: \G \to \R_+$ is a length function if and only if the mappings $S_{\psi,t}(\lambda_\G(g)) = \exp(-t\psi(g)) \lambda_\G(g)$ extend to a semigroup of unital completely positive maps on $\V$. Let us collect these results.

\begin{lemma} \label{Lemmapsi}
If $\psi: \G \to \R_+$ is a length function$\, :$
\begin{itemize}
\item The form
\begin{eqnarray*}
K_\psi(g,h) & = & \frac{\psi(g) + \psi(h) - \psi(g^{-1}h)}{2},
\end{eqnarray*}
defines a positive matrix on $\G \times \G$ and leads to $$\Big\langle \summ_g a_g \delta_g, \summ_h b_h
\delta_h \Big\rangle_{\psi} = \summ_{g,h} a_g K_\psi(g,h) b_h$$ on the group algebra $\R[\G]$ of finitely supported real functions on $\G$.

\vskip5pt

\item Let $\H$ be the Hilbert space completion of $$(\R[\G]/N_\psi, \langle \cdot \hskip1pt, \cdot \rangle_{\psi}) \quad \mbox{with} \quad N_\psi = \mbox{null space of} \ \langle \cdot \hskip1pt , \cdot \rangle_{\psi}.$$ If we consider the mapping $b: g \in \G \mapsto \delta_g + N_\psi \in \H$
\begin{eqnarray*}
\gamma_g \Big( \sum_{h \in \G} a_h b(h) \Big) & = & \sum_{h \in \G} a_h \big( b(gh) - b(g) \big)
\end{eqnarray*}
determines an isometric action $\gamma: \G \to O(\H)$ of $\G$ on $\H$.

\vskip5pt

\item The resulting triple $(\H, \gamma, b)$ determines a cocycle of the group $\G$.
\end{itemize}
\end{lemma}

The previous lemma allows the reader to consider a pseudo-metric on the discrete group $\G$ in terms of the length function $\psi$. Indeed, a short calculation leads to the crucial identity $\psi(g^{-1}h) = \langle b(g) - b(h), b(g) - b(h) \rangle_{\psi} = \| b(g) - b(h) \|_{\H}^2.$ In particular we find that $\mathrm{dist}(g,h) = \sqrt{\psi(g^{-1}h)} = \|b(g) - b(h) \|_{\H}$ defines a pseudo-metric on $\G$, which becomes a metric when the $b$ is injective.

\subsection{Littlewood-Paley estimates for length functions}

Consider a family of operators $\Sigma = (\sigma_k)_{k \ge 1}$ acting on some Hilbert space. Then, the row and column square functions associated to $\Sigma$ are respectively defined by $$\mathcal{S}_r(\Sigma) = \Big( \summ_k \sigma_k \sigma_k^* \Big)^\frac12 \qquad \mbox{and} \qquad \mathcal{S}_c(\Sigma) = \Big( \summ_k \sigma_k^* \sigma_k \Big)^\frac12.$$ Given a discrete group $\G$ and a family $\Sigma = (\sigma_k)_{k \ge 1}$ in $L_p(\widehat{\mathbf{G}})$, consider the norms $\|\Sigma\|_{L_p(\widehat{\mathbf{G}}; \ell_2^r)} = \| \mathcal{S}_r(\Sigma) \|_p$ and $\|\Sigma\|_{L_p(\widehat{\mathbf{G}}; \ell_2^c)} = \| \mathcal{S}_c(\Sigma) \|_p$. Both clearly coincide over commutative algebras. In general, certain combination is needed to obtain the noncommutative forms of classical results such as Khintchine, Burholder-Gundy or Littlewood-Paley type inequalities. It is now well-known that the right combination arises as follows $$L_p(\widehat{\mathbf{G}}; \ell_{rc}^2) = \begin{cases} L_p(\widehat{\mathbf{G}}; \ell_2^r) + L_p(\widehat{\mathbf{G}}; \ell_2^c) & \mbox{if} \ 1 \le p \le 2, \\ L_p(\widehat{\mathbf{G}}; \ell_2^r) \cap \hskip1pt L_p(\widehat{\mathbf{G}}; \ell_2^c) & \mbox{if} \ 2 \le p \le \infty. \end{cases}$$ In other words, we have $$\|\Sigma\|_{L_p(\widehat{\mathbf{G}},\ell_{rc}^2)} = \begin{cases} \displaystyle \inf_{\Sigma = \Phi + \Psi} \|\mathcal{S}_r(\Phi)\|_p + \|\mathcal{S}_c(\Psi)\|_p & \mbox{if} \ 1 \le p \le 2, \\ \displaystyle \ \, \max \ \big\{ \|\mathcal{S}_r(\Sigma)\|_p, \|\mathcal{S}_c(\Sigma)\|_p \big\} & \mbox{if} \ 2 \le p \le \infty. \end{cases}$$ The following result can be found in \cite{JMP1}.

\begin{lemma} \label{Lema3}
Let $\G$ be a discrete group equipped with a length function $\psi: \G \to \R_+$ and assume that $\dim \H = n < \infty$ for the  cocycle Hilbert space. Let $k_n = [\frac{n}{2}] + 1$ and consider a family $h_m \in \mathcal{C}^{k_n}(\R_+ \setminus \{0\})$ satisfying
\begin{itemize}
\item $\summ_m |h_m(\zeta)|^2 = 1$,

\vskip3pt

\item $\summ_m \big| \frac{d^j}{d\zeta} h_m(\zeta) \big|^2  \le c_n |\zeta|^{-2j}$ for $j \le [\frac{n}{2}] + 1$.
\end{itemize}
Then, the following holds for $f \in L_p(\widehat{\mathbf{G}})$ and $1 < p < \infty$ $$\|f\|_p \, \sim_{c_n} \, \Big\| \summ_m f_m \otimes e_m \Big\|_{L_p(\widehat{\mathbf{G}};\ell_{rc}^2)},$$ with the $\psi$-smooth Littlewood-Paley decomposition $f_m = \sum_g h_m(\psi(g)) \widehat{f}(g) \lambda_\G(g)$.
\end{lemma}

\subsection{Length functions on $\G$ adapted to $\gamma$}
\label{PLength-gamma}

It may be illustrative for non-experts to show how to construct natural length functions for those discrete groups which admit finite-dimensional orthogonal representations. This will be used below in the proof of Theorem B. Take $$\psi_{O(n)}(\mathrm{A}) \, = \, \big\| \mathrm{A} - \mathrm{I} \big\|_{\mathrm{HS}}^2 \, = \, \summ_{jk} \big| \mathrm{A}_{jk} - \delta_{jk} \big|^2.$$ It is not difficult to check directly that $\psi_{O(n)}$ is a length function in $O(n)$, but it is perhaps easier to note that $\psi_{O(n)}(\mathrm{A}) = \|b(\mathrm{A})\|_{\H}^2$ for the cocycle $(\H,\gamma,b)$ which is determined by the $n \times n$ matrices with the Hilbert-Schmidt norm, the action $\gamma_\mathrm{A}(\mathrm{B}) = \mathrm{AB}$ and the cocycle map $\mathrm{A} \mapsto \mathrm{A-I}$. Consider now a discrete group equipped with a orthogonal representation $\gamma: \G \to O(n)$. Using that $\psi_{O(n)}$ is a length function, we may define $\psi_\gamma: \G \to \R_+$ as follows $$\psi_\gamma(g) \, = \, \psi_{O(n)}(\gamma_g) \, = \, \big\| \gamma_g - \mathrm{I} \big\|_{\mathrm{HS}}^2.$$ It is now clear that $\psi_\gamma$ defines a length function on the group $\G$ for any orthogonal representation $\gamma$. Alternatively, given any non-zero $\xi_0 \in \R^n$ we may also construct the length functions $$\psi_{\gamma,\xi_0}(g) \, = \, \big\langle \gamma_g(\xi_0) - \xi_0, \gamma_g(\xi_0) - \xi_0 \big\rangle_{\R^n}.$$ Both choices of length functions correspond to inner cocycles ---$b(g) = \gamma_g (\eta) - \eta$ for some $\eta \in \H$--- which are quotiented out in the formation of the corresponding cohomology group, so the reader could object that our length functions are singular in the sense of cohomology theory. This was already justified in \cite{JMP1}, where inner cocycles turned out to be the most striking ones looking for pathological Fourier multipliers, even in $\R^n$. On the other hand, discrete groups satisfying Kazhdan's property $(\mathrm{T})$ only admit inner cocycles. Thus, the information encoded by our length functions goes beyond the cohomology group, which is trivial for this class of groups. Bounded, integer valued lengths also arise regarding $O(n)$ as a Coxeter group, counting the number of reflections in which $\gamma_g$ decomposes.

\section{Twisted Hilbert transforms vs Kakeya sets}

In this section we prove Theorem A. Most of our efforts are devoted to proving the hardest implication ii) $\Rightarrow$ iii). First, we use Littlewood-Paley estimates in a group amplification of $\G$ to obtain a square function inequality in $L_p(\R^n_{\mathrm{bohr}})$, provided $H_u \rtimes_\gamma id_\G$ is $L_p$-bounded. The group amplification is essential to provide enough room to disprove such an inequality for infinite orbits. Second, we will decompactify such an inequality adapting transference techniques in conjunction with ergodic type arguments, which yields an inequality in $L_p(\R^n)$. Third, we show that infinite orbits admit Kakeya shadows and disprove the latter Euclidean inequality. The proof of i) $\Rightarrow$ ii) uses a slight variation of de Leeuw's compactification argument which adapts to the cross product setting. Finally, the $L_p$-boundedness for finite orbits is clear and we shall prove stronger endpoint estimates.

\subsection{Littlewood-Paley methods}

The following result plays a role similar to Meyer's lemma in Fefferman's solution of the disc conjecture \cite{Fef}. By duality in Theorem A, it suffices to consider the case $2 < p < \infty$. However, the inequality in the following lemma is not self-dual. One can formulate a (more intricate) version for $1 < p < 2$ which will be easily guessed by the reader after the proof.

\begin{lemma} \label{LPReduction1}
If $2 < p < \infty$ and $H_u \rtimes_\gamma id_\G: L_p(\widehat{\mathbf{\Gamma}_\mathrm{disc}}) \to L_p(\widehat{\mathbf{\Gamma}_\mathrm{disc}})$, then $$\Big\| \Big( \sum_{j=1}^\infty |H_{\gamma_{g_j}(u)}(f_{g_j})|^2 \Big)^\frac12 \Big\|_p \, \le \, c_p \max \left\{ \Big\| \Big( \sum_{j=1}^\infty |f_{g_j}|^2 \Big)^\frac12 \Big\|_p, \Big\| \Big( \sum_{j=1}^\infty |\gamma_{g_j^{-1}}f_{g_j}|^2 \Big)^\frac12 \Big\|_p \right\}$$ for any sequence $g_1, g_2, \ldots$ in $\G$ and any family of functions $f_{g_1}, f_{g_2}, \ldots$ in $L_p(\R_{\mathrm{bohr}}^n)$.
\end{lemma}

\dem Let $\mathrm{H} = \G \times \Z$ and set $$\rho: (g,k) \in \mathrm{H} \mapsto \gamma_g \in O(n).$$ Clearly, $\rho$ defines an orthogonal representation of $\mathrm{H}$ which yields
\begin{eqnarray*}
\R_\mathrm{disc}^n \rtimes_\rho \mathrm{H} & \simeq & (\R_\mathrm{disc}^n \rtimes_\gamma \G) \times \Z, \\ \mathcal{L}(\R_\mathrm{disc}^n \rtimes_\rho \mathrm{H}) & \simeq & \mathcal{L}(\R_\mathrm{disc}^n \rtimes_\gamma \G) \bar\otimes \mathcal{L}(\Z).
\end{eqnarray*}
This group isomorphism $\xi \rtimes_\rho (g,k) \mapsto (\xi \rtimes_\gamma g, k)$ is additionally a homeomorphism since both sides are equipped with the discrete topology. The map $H_u \rtimes_\rho id_\mathrm{H}$ factorizes as $(H_u \rtimes_\gamma id_\G) \otimes id_\Z$, so that $L_p$-boundedness means $$\int_\mathbb{T} \big\| (H_u \rtimes_\gamma id_\G)f(x) \big\|_p^p \, dx \, \le \, c_p \int_\mathbb{T} \|f(x)\|_p^p \, dx$$ for any $f \in L_p(\mathbb{T}; L_p(\widehat{\mathbf{\Gamma}_{\mathrm{disc}}}))$, which clearly holds by hypothesis. Once we know $H_u \rtimes_\rho id_\mathrm{H}$ is $L_p$-bounded, consider the length function on $\R_\mathrm{disc}^n \rtimes_\rho \mathrm{H}$ given by $\psi (\xi \rtimes_\rho (g,k)) = |k|^2$ for all $(g,k,\xi) \in \G \times \Z \times \R^n$. Recall that $\psi$ gives rise to the trivial cocycle $b(\xi \rtimes_\rho (g,k)) = k$. Thus, we may apply Lemma \ref{Lema3} with $n = 1$. If we pick a small $\delta > 0$ and any nonnegative radially decreasing Schwartz function $\phi: \R \to \R_+$ which takes the value $1$ for $|\zeta| \le \frac{1 + \delta}{2}$ and vanishes for $|\zeta| \ge 1-\delta$, the family of functions $$h_m(\zeta) = \Big( \phi \big( \frac{2^{-m}\zeta}{2} \big) - \phi(2^{-m} \zeta) \Big)^\frac12 \quad \mbox{for} \quad m \in \Z$$ trivially satisfy the hypotheses of the lemma. Now, given elements $g_1, g_2, \ldots$ in $\G$ and functions $f_{g_1}, f_{g_2}, \ldots$ in $L_p(\R^n_{\mathrm{disc}})$, we set $f = \summ_j \phi_{g_j} \rtimes_\rho \lambda_\mathrm{H}(g_j^{-1},2^j)$ with $\phi_{g_j} = \gamma_{g_j^{-1}} f_{g_j}$ and recall the identities
\begin{eqnarray*}
\summ_j \big( \phi_{g_j} \rtimes_\rho \lambda_\mathrm{H}(g_j^{-1},2^j) \big) \big( \phi_{g_j} \rtimes_\rho \lambda_\mathrm{H}(g_j^{-1},2^j) \big)^* & = & \Big(\summ_j | \phi_{g_j} |^2 \Big) \rtimes_\rho \mathbf{1}_{\mathcal{L}(\mathrm{H})}, \\ \summ_j \big( \phi_{g_j} \rtimes_\rho \lambda_\mathrm{H}(g_j^{-1},2^j) \big)^* \big( \phi_{g_j} \rtimes_\rho \lambda_\mathrm{H}(g_j^{-1},2^j) \big) & = & \Big(\summ_j |\gamma_{g_j} \phi_{g_j}|^2 \Big) \rtimes_\rho \mathbf{1}_{\mathcal{L}(\mathrm{H})}.
\end{eqnarray*}
According to Lemma \ref{Lema3}, we find $$\hskip-4pt \|f\|_{L_p(\mathcal{L}(\R^n_{\mathrm{disc}} \rtimes_\rho \mathrm{H}))} \, \sim \, \max \Big\{ \Big\| \Big(\summ_j |f_{g_j}|^2 \Big)^\frac12 \Big\|_p, \Big\| \Big(\summ_j |\gamma_{g_j^{-1}}f_{g_j}|^2 \Big)^\frac12 \Big\|_p \Big\}.$$ Using the same norm equivalence for $(H_u \rtimes_\rho id_\mathrm{H})f = H_u[f]$, we also get
$$\|H_u[f]\|_{L_p(\mathcal{L}(\R^n_{\mathrm{disc}} \rtimes_\rho \mathrm{H}))} \, \sim \, \max \Big\{ \Big\| \Big(\summ_j |H_u \phi_{g_j}|^2 \Big)^\frac12 \Big\|_p, \Big\| \Big(\summ_j |\gamma_{g_j}H_u \phi_{g_j}|^2 \Big)^\frac12 \Big\|_p \Big\}.$$ By de Leeuw's theorem, $H_u$ is bounded on $L_p(\R^n_{\mathrm{disc}})$. Thus, the first term on the right hand side for $H_u[f]$ is dominated by the second term on the right hand side for $f$. Therefore, since $H_u \rtimes_\rho id_\mathrm{H}$ is $L_p$-bounded, the second term on the right hand side for $H_u[f]$ must be dominated by the maximum associated to $f$. Now, using the identity $\gamma_g H_u \phi = H_{\gamma_g(u)}\gamma_g \phi$ we recover the desired inequalities. \fin

\begin{remark} \label{Bourgain}
\emph{The full strength of Lemma \ref{Lema3} is not necessary to prove Lemma \ref{LPReduction1}. One can also combine Bourgain's extension of Littlewood-Paley estimates for UMD Banach spaces \cite{Bo} together with Lust-Piquard/Pisier's noncommutative Khintchine inequalities \cite{L,LP}. However, we will require Lemma \ref{Lema3} later to prove Theorem B, and  we find our approach more intrinsic, which could help for future generalizations.}
\end{remark}

\subsection{A partial decompactification} \label{Decompact}

The inequality arising from Lemma \ref{LPReduction1} is stated for $L_p$-functions on the Bohr compactification of $\R^n$, with respect to the corresponding normalized Haar measure. In this paragraph we are interested in the form that such an inequality takes in $L_p(\R^n)$ with the Lebesgue measure. Given $\mathrm{M} > 0$ and $f: \R^n \to \C$ supported by $[0, \mathrm{M}] \times \cdots \times [0, \mathrm{M}]$, we shall write $\pi_\mathrm{M}f$ for its $\mathrm{M}$-periodization in the axes directions $\pi_\mathrm{M}f(x) = \sum_{k \in \Z^n} f (x - \mathrm{M} k)$. We also consider $\mathrm{M}$-periodizations along $g$-lattices $$\pi_\mathrm{M}^gf(x) = \sum_{k \in \Z^n} f \big( x - \mathrm{M} \gamma_{g^{-1}}(k) \big).$$ Let us also recall the $L_p$ norm for almost periodic functions in $\R^n$ $$\left\bracevert \phi \right\bracevert_p \, = \, \Big( \lim_{\Lambda \to \infty} \frac{1}{(2\Lambda)^n} \int_{[-\Lambda,\Lambda]^n} |\phi(x)|^p \, dx \Big)^{\frac{1}{p}}.$$

\begin{lemma} \label{Decompactification}
If $2 < p < \infty$ and $H_u \rtimes_\gamma id_\G: L_p(\widehat{\mathbf{\Gamma}_\mathrm{disc}}) \to L_p(\widehat{\mathbf{\Gamma}_\mathrm{disc}})$, then
\begin{eqnarray*}
\lefteqn{\hskip-10pt \Big\| \Big( \sum_{j=1}^\infty |H_{\gamma_{g_j}(u)}(f_{g_j})|^2 \Big)^\frac12 \Big\|_p} \\ & \le & c_p \max \left\{ \Big\| \Big( \sum_{j=1}^\infty |f_{g_j}|^2 \Big)^\frac12 \Big\|_p, \lim_{\mathrm{M} \to \infty} \mathrm{M}^{\frac{n}{p}} \Big\bracevert \Big( \sum_{j=1}^\infty |\pi_\mathrm{M}^{g_j} \gamma_{g_j^{-1}} f_{g_j}|^2 \Big)^\frac12 \Big\bracevert_p \right\}
\end{eqnarray*}
for any sequence $g_1, g_2, \ldots$ in $\G$ and any family of functions $f_{g_1}, f_{g_2}, \ldots$ in $L_p(\R^n)$.
\end{lemma}

\dem Our argument adapts de Leeuw's decompactification argument in \cite{dL}. By density, it suffices to prove such an inequality for a finite family $f_{g_1}, f_{g_2}, \ldots, f_{g_\mathrm{N}}$ of compactly supported Schwartz functions. If $S_p$ denotes the Schatten $p$-class, the left hand side is the norm in $L_p(\R^n; S_p)$ of $\summ_j H_{\gamma_j(u)}(f_{g_j}) \otimes e_{j1}$. In particular, trace duality provides us with a (matrix-valued) compactly supported Schwartz function $h = \sum_j h_j \otimes e_{j1}$ in the unit ball of $L_q(\R^n; S_q)$ for $\frac1p + \frac1q = 1$, such that
\begin{eqnarray*}
\Big\| \Big( \sum_{j=1}^\mathrm{N} |H_{\gamma_{g_j}(u)}(f_{g_j})|^2 \Big)^\frac12 \Big\|_p & \sim & \sum_{j=1}^\mathrm{N} \int_{\R^n} H_{\gamma_{g_j}(u)}(f_{g_j})(x) \overline{h_j(x)} \, dx \\ & = & \sum_{j=1}^\mathrm{N} \int_{\R^n} \widehat{H_{\gamma_{g_j}(u)}(f_{g_j})}(\xi) \overline{\widehat{h_j}(-\xi)} \, d\xi \\ & = & \sum_{j=1}^\mathrm{N} \lim_{\mathrm{M} \to \infty} \frac{1}{\mathrm{M}^n} \sum_{k \in \frac{1}{\mathrm{M}} \Z^n} - i \mathrm{sgn} \langle \gamma_{g_j}(u), k \rangle \widehat{f_{g_j}}(k) \overline{\widehat{h_j}(-k)}.
\end{eqnarray*}
If $\mathrm{supp} f \subset [-\frac{\mathrm{M}}{2}, \frac{\mathrm{M}}{2}]^n$, we recall the identity $$\widehat{f}(k) = \mathrm{M}^n \widehat{\pi_\mathrm{M}f}(k),$$ with the first Fourier transform calculated in $\R^n$ and the second in $\mathrm{M} \mathbb{T}^n$. This gives $$\Big\| \Big( \sum_{j=1}^\mathrm{N} |H_{\gamma_{g_j}(u)}(f_{g_j})|^2 \Big)^\frac12 \Big\|_p \, \sim \, - i \lim_{\mathrm{M} \to \infty} \mathrm{M}^n \sum_{j,k} \mathrm{sgn} \langle \gamma_{g_j}(u), k \rangle \widehat{\pi_\mathrm{M}f_{g_j}}(k) \overline{\widehat{\pi_\mathrm{M}h_j}(-k)}.$$ Let us write $\exp_\xi: \R^n \to \mathbb{T}$ for the usual characters of $\R^n$ and $\mbox{b-exp}_\xi: \R^n_{\mathrm{bohr}} \to \mathbb{T}$ for the characters of the Bohr compactification. Both families are indexed by the same set, but the latter are defined in a larger group. Define $$\phi_{\mathrm{M},j} = \sum_{\xi \in \frac{1}{\mathrm{M}} \Z^n} \widehat{\pi_\mathrm{M}f_{g_j}}(\xi) \mbox{b-exp}_\xi \quad \mbox{and} \quad \varphi_{\mathrm{M},j} = \sum_{\xi \in \frac{1}{\mathrm{M}} \Z^n} \widehat{\pi_\mathrm{M} \hskip1pt h_j}(\xi) \hskip1pt \mbox{b-exp}_\xi.$$
This yields
\begin{eqnarray*}
- i \summ_k \mathrm{sgn} \langle \gamma_{g_j}(u), k \rangle \widehat{\pi_\mathrm{M}f_{g_j}}(k) \overline{\widehat{\pi_\mathrm{M}h_j}(-k)} & = & \summ_\xi \widehat{H_{\gamma_{g_j}(u)} (\phi_{\mathrm{M},j})}(\xi) \overline{\widehat{\varphi_{\mathrm{M},j}}(-\xi)} \\ & = & \int H_{\gamma_{g_j}(u)} (\phi_{\mathrm{M},j})(x) \overline{\varphi_{\mathrm{M},j}(x)} \, d\mu(x),
\end{eqnarray*}
where now the Hilbert transforms must be understood as operators on the Bohr compactification and $\mu$ denotes the corresponding normalized Haar measure. On the other hand, trace duality in $L_2(\R^n_{\mathrm{bohr}}; S_2)$ gives
\begin{eqnarray*}
\lefteqn{\Big| \summ_j \int H_{\gamma_{g_j}(u)} (\phi_{\mathrm{M},j})(x) \overline{\varphi_{\mathrm{M},j}(x)} \, d\mu(x) \Big|} \\ & \le & \Big\| \summ_j H_{\gamma_{g_j}(u)}(\phi_{\mathrm{M},j}) \otimes e_{j1} \Big\|_{L_p(\R^n_{\mathrm{bohr}};S_p)} \Big\| \summ_j \varphi_{\mathrm{M},j} \otimes e_{j1} \Big\|_{L_q(\R^n_{\mathrm{bohr}};S_q)} \ = \ \mathrm{AB}.
\end{eqnarray*}
According to Lemma \ref{LPReduction1}, we know that $$\mathrm{A} \, \le \, c_p \max \left\{ \Big\| \summ_j \phi_{\mathrm{M},j} \otimes e_{j1} \Big\|_{L_p(\R^n_{\mathrm{bohr}};S_p)}, \Big\| \summ_j\gamma_{g_j^{-1}}\phi_{\mathrm{M},j} \otimes e_{j1} \Big\|_{L_p(\R^n_{\mathrm{bohr}};S_p)} \right\}.$$ Let us now recall the identity for trigonometric polynomials
$$\Big\| \summ_\xi \widehat{f}(\xi) \mbox{b-exp}_\xi \Big\|_{L_p(\R^n_{\mathrm{bohr}})} \ = \ \Big\bracevert \summ_\xi \widehat{f}(\xi) \exp_\xi \Big\bracevert_p$$ which relates Haar integration and mean values. Approximating by trigonometric polynomials
and noticing again the identity between $L_p$-norms of column matrices and $L_p$-norms of square functions, we find
\begin{eqnarray*}
\Big\| \summ_j \phi_{\mathrm{M},j} \otimes e_{j1} \Big\|_{L_p(\R^n_{\mathrm{bohr}};S_p)} & = & \Big\bracevert \Big( \summ_j |\pi_\mathrm{M} f_{g_j}|^2 \Big)^\frac12 \Big\bracevert_p \\ & = & \hskip2pt \Big\| \Big( \summ_j |\pi_\mathrm{M} f_{g_j}|^2 \Big)^\frac12 \Big\|_{L_p(\mathrm{M} \mathbb{T}^n)} \\ & = & \frac{1}{\mathrm{M}^{\frac{n}{p}}} \Big\| \Big( \summ_j |f_{g_j}|^2 \Big)^\frac12 \Big\|_{L_p(\mathbb{R}^n)}.
\end{eqnarray*}
The second identity follows by $\mathrm{M}$-periodicity and the third one is valid for $\mathrm{M}$ large enough, since the $f_{g_j}$'s are compactly supported. The same argument and the fact that $\sum_j h_j \otimes e_{j1}$ is in the unit ball of $L_q(\R^n;S_q)$ gives $\mathrm{B} = \mathrm{M}^{- \frac{n}{q}}$. Finally, since
\begin{eqnarray*}
\gamma_{g_j^{-1}} \phi_{\mathrm{M},j} & = & \sum_{\mathrm{M} \xi \in \Z^n} \widehat{\pi_\mathrm{M}f_{g_j}}(\xi) \mbox{b-exp}_{\gamma_{g_j^{-1}}\xi} \ = \ \sum_{\mathrm{M} \xi \in \Z^n} \frac{\widehat{f_{g_j}}(\xi)}{\mathrm{M}^n} \mbox{b-exp}_{\gamma_{g_j^{-1}}\xi} \\ & = & \sum_{\mathrm{M} \xi \in \gamma_{g_j^{-1}}(\Z^n)} \frac{\widehat{\gamma_{g_j^{-1}} f_{g_j}}(\xi)}{\mathrm{M}^n} \mbox{b-exp}_{\xi}  \ = \ \sum_{\mathrm{M} \xi \in \gamma_{g_j^{-1}}(\Z^n)} \widehat{\pi_\mathrm{M}^{g_j}\gamma_{g_j^{-1}} f_{g_j}}(\xi) \mbox{b-exp}_{\xi},
\end{eqnarray*}
we obtain the following identity as before $$\Big\| \summ_j\gamma_{g_j^{-1}}\phi_{\mathrm{M},j} \otimes e_{j1} \Big\|_{L_p(\R^n_{\mathrm{bohr}};S_p)} = \Big\bracevert \Big( \summ_j |\pi_\mathrm{M}^{g_j} \gamma_{g_j^{-1}} f_{g_j}|^2 \Big)^\frac12 \Big\bracevert_p.$$ The assertion now follows by combining the estimates obtained so far. \fin

\subsection{Distribution of points in lattice intersections} \label{PLattice}

The decompactification of the last term in Lemma \ref{Decompactification} requires a more careful analysis. Let us consider the  two-dimensional lattice $\mathrm{M} \rho_{\alpha}(\Z^2)$ given by an $\alpha$-rotation of $\mathrm{M} \Z^2$. We set $\pi_\mathrm{M}^\alpha$ for the corresponding periodization operator and $$\langle \beta \rangle = \begin{cases} \infty & \mbox{if} \ \beta \in \R \backslash \mathbb{Q}, \\ \sqrt{p^2 + q^2} & \mbox{if} \ \beta = p/q \in \mathbb{Q}, \end{cases}$$ where $p/q$ is written in reduced form, so that $p,q$ are relatively prime.

\begin{lemma} \label{ErgodicLemma}
Given $z_0 \in \R^2$ and $\delta > 0$
$$\lim_{\mathrm{M} \to \infty} \lim_{\Lambda \to \infty} \frac{\mathrm{M}^2}{4\Lambda^2} \Big| [-\Lambda,\Lambda]^2 \cap \pi_{\mathrm{M}} (\mathrm{B}_\delta(z_0)) \cap \pi_\mathrm{M}^\alpha (\mathrm{B}_\delta(z_0)) \Big| \, \lesssim \, \frac{\delta^2}{\langle \tan \alpha \rangle}.$$
\end{lemma}

\dem Classical ergodic theory \cite{Do} gives that $$\mathrm{Orb}_\mathrm{M}(\alpha) \, = \, \Big\{ \big( x, y) \, \mbox{mod} \, \mathrm{M} \times \mathrm{M} \ \big| \ (x,y) \in \pi_\mathrm{M}^\alpha(\{z_0\}) \Big\}$$ is injective, dense and uniformly distributed in the unit cell $[0,\mathrm{M}] \times [0,\mathrm{M}]$ when $\tan \alpha \in \R \backslash \mathbb{Q}$. In particular, the following estimate holds for irrational slopes and sufficiently large $\mathrm{M}$ $$\lim_{\Lambda \to \infty} \frac{\mathrm{M}^2}{4\Lambda^2} \Big| [-\Lambda,\Lambda]^2 \cap \pi_{\mathrm{M}} (\mathrm{B}_\delta(z_0)) \cap \pi_\mathrm{M}^\alpha (\mathrm{B}_\delta(z_0)) \Big| \ \lesssim \ \frac{\delta^4}{\mathrm{M}^2}.$$ Indeed, let $\mathcal{Q}_\Lambda$ be a covering of $[-\Lambda,\Lambda]^2$ by disjoint adjacent $\mathrm{M} \times \mathrm{M}$ cubes with sides parallel to the axes. Note that $|\mathcal{Q}_\Lambda| \sim 4\Lambda^2 / \mathrm{M}^2$. On the other hand, if we pick a random point $\sigma$ in $\mathrm{M} \rho_\alpha(\Z^2)$, the probability that $\mathrm{B}_\delta(z_0) \cap (\mathrm{B}_\delta(z_0) + \sigma \ \mathrm{mod} \ \mathrm{M} \times \mathrm{M})$ is not empty is the same that $z_0 + \sigma \ \mathrm{mod} \ \mathrm{M} \times \mathrm{M}$ belongs to the ball $\mathrm{B}_{2\delta}(z_0)$, which in turn is comparable to $\frac{1}{\mathrm{M}^2} |\mathrm{B}_{2\delta}(z_0)|$ since the $\pi_\mathrm{M}^\alpha$-orbit of $z_0$ is uniformly distributed in the unit cell. In conclusion, we may rewrite the left-hand side as follows $$\lim_{\Lambda \to \infty} \frac{1}{|\mathcal{Q}_\Lambda|} \sum_{Q \in \mathcal{Q}_\Lambda} \Big| Q \cap \pi_\mathrm{M}(\mathrm{B}_\delta(z_0)) \cap \pi_\mathrm{M}^\alpha(\mathrm{B}_\delta(z_0)) \Big|.$$ Translating everything to the unit cell, we get
\begin{eqnarray*}
\lefteqn{\hskip-10pt \lim_{\Lambda \to \infty} \frac{1}{|\mathcal{Q}_\Lambda|} \sum_{Q \in \mathcal{Q}_\Lambda} \Big| \mathrm{B}_\delta(z_0) \cap \big( \pi_\mathrm{M}^\alpha(\mathrm{B}_\delta(z_0)) \cap Q \ \mathrm{mod} \ \mathrm{M} \times \mathrm{M} \big) \Big|} \\ & \lesssim & |\mathrm{B}_\delta(z_0)| \, \mathrm{Prob} \Big\{ \mathrm{B}_\delta(z_0) \cap (\mathrm{B}_\delta(z_0) + \sigma \ \mathrm{mod} \ \mathrm{M} \times \mathrm{M}) \neq \emptyset \Big\} \ \sim \ \frac{\delta^4}{\mathrm{M}^2}.
\end{eqnarray*}
Taking limits in $\mathrm{M}$, we conclude for $\tan \alpha \in \R \backslash \mathbb{Q}$. If $\tan \alpha = p/q$ is rational, set $$e_1 = \frac{(p,q)}{\sqrt{p^2+q^2}} \quad \mbox{and} \quad e_2 = \frac{(-q,p)}{\sqrt{p^2+q^2}}.$$ Of course, we choose $p,q$ relatively prime. Now we may write $$\mathrm{Orb}_\mathrm{M}(\alpha) \, = \, \Big\{ \big( z_0 + j \mathrm{M}e_1 + k \mathrm{M} e_2 \big) \, \mbox{mod} \, \mathrm{M} \times \mathrm{M} \ \big| \ j,k \in \Z \Big\}.$$ Therefore, the suborbits generated by $\mathrm{M}e_1$ and $\mathrm{M}e_2$ span $\mathrm{Orb}_\mathrm{M}(\alpha)$. Note that the suborbit generated by $\mathrm{M} e_j$ lives inside a union of segments $\Sigma_k$ in the unit cell which point in the direction of $e_j$ and such that $$\summ_k |\Sigma_k| = \mathrm{M} \sqrt{p^2+q^2}.$$ When this square root is irrational, $\mathrm{Orb}_\mathrm{M}(\alpha)$ is again injective, dense and uniformly distributed, and the argument for irrational slopes still applies. When it is rational the orbit $\mathrm{Orb}_\mathrm{M}(\alpha)$ is finite and its cardinality $J_\alpha = |\mathrm{Orb}_\mathrm{M}(\alpha)|$ is independent of $\mathrm{M}$. In particular, we set $$\Delta_\mathrm{M}(\alpha) = \inf \Big\{ \mathrm{dist}(\mathrm{A},\mathrm{B}) \, \big| \ \mathrm{A}, \mathrm{B} \in \mathrm{Orb}_\mathrm{M}(\alpha) \Big\} > 0.$$ Note that $\Delta_\mathrm{M}(\alpha) \sim \frac{\mathrm{M}}{\mathrm{M}_0} \Delta_{\mathrm{M}_0}(\alpha)$ for $\mathrm{M} \ge \mathrm{M}_0$ large enough, so that $$\lim_{\mathrm{M} \to \infty} \lim_{\Lambda \to \infty} \frac{\mathrm{M}^2}{4\Lambda^2} \Big| [-\Lambda,\Lambda]^2 \cap \pi_{\mathrm{M}} (\mathrm{B}_\delta(z_0)) \cap \pi_\mathrm{M}^\alpha (\mathrm{B}_\delta(z_0)) \Big| \lesssim \frac{|\mathrm{B}_\delta(z_0)|}{J_\alpha} \sim \frac{\delta^2}{J_\alpha}$$ since we eventually find $\Delta_\mathrm{M}(\alpha)> 2 \delta$. However, the suborbit generated by $\mathrm{M}e_1$ contains at least $\sqrt{p^2+q^2}$ nonequivalent points, more if $\sqrt{p^2+q^2} \in \mathbb{Q} \backslash \Z$. In particular, we always have $J_\alpha \ge \langle \tan \alpha \rangle$ and the proof is complete. \fin

For the last term in Lemma \ref{Decompactification}, we will consider a finite family of pairwise commuting $g_j \in \G$, with $1 \le j \le \mathrm{N}$. In particular, the orthogonal maps $\gamma_{g_j}$ will admit a simultaneous diagonal form. This means that we may find a direct sum decomposition $$\R^n = \Delta_{+ 1} \oplus \Delta_{- 1} \oplus \Theta_1 \oplus \Theta_2 \oplus \cdots \oplus \Theta_\ell,$$ where $\Delta_{\pm 1}$ is the direct sum of the eigenspaces with eigenvalues $\pm 1$ and the $\Theta_k$'s are $2$-dimensional eigenspaces where the $\gamma_{g_j}$'s act by rotations. On the other hand, our functions $f_{g_j}$ will be characteristic functions of prisms $$A_j = R_j \times [-\lambda,\lambda]^{n-2}$$ with $R_j$ disjoint rectangles living in $\Theta_1$ and certain $\lambda > 0$ so that $\summ_j |A_j| \ge 1$. Let us write $\rho_{g_j}$ and $\rho_{g_j,\perp}$ for the restriction of $\gamma_{g_j}$ to $\Theta_1$ and its orthocomplement respectively. We write $\alpha_{jk}$ for the rotation angle of the map $\rho_{g_jg_k^{-1}}$.

\begin{lemma} \label{Ergodic2}
Given $p_0 \ge 2$, we have
$$\lim_{\mathrm{M} \to \infty} \mathrm{M}^{\frac{n}{p_0}} \Big\bracevert \Big( \sum_{j=1}^\mathrm{N} |\pi_\mathrm{M}^{g_j} \gamma_{g_j^{-1}} \chi_{A_j}|^2 \Big)^\frac12 \Big\bracevert_{p_0} \, \sim \, \Big| \bigcup_{j=1}^\mathrm{N} A_j \Big|^{\frac{1}{p_0}}.$$ provided that $\big\langle \tan \alpha_{jk} \big\rangle \displaystyle \ge \mathrm{N}^{p_0+1} \lambda^{n-2} \mathrm{diam}^2 \Big( \bigcup_{j=1}^\mathrm{N} \rho_j^{-1}(R_j) \Big)$ for any $1 \le j \neq k \le \mathrm{N}$.
\end{lemma}

\dem If $\mathrm{M}$ is large enough $\pi_\mathrm{M}^{g_j} \gamma_{g_j}^{-1} \chi_{A_j} = \chi_{\mathrm{B}_{\mathrm{M},j}}$, where $\mathrm{B}_{\mathrm{M},j}$ is the periodization of $\gamma_{g_j}^{-1} (\mathrm{A}_j)$ along $\mathrm{M} \gamma_{g_j}^{-1}(\Z^n)$. If the periodized sets $\mathrm{B}_{\mathrm{M},j}$ were pairwise disjoint for all $\mathrm{M} \ge \mathrm{M}_0$, the assertion would trivially follow by direct computation. In the presence of overlapping $$\mathrm{M}^{\frac{n}{p_0}} \Big\bracevert \Big( \sum_{j=1}^\mathrm{N} |\pi_\mathrm{M}^{g_j} \gamma_{g_j^{-1}} \chi_{A_j}|^2 \Big)^\frac12 \Big\bracevert_{p_0} \, \ge \, \Big| \bigcup_{j=1}^\mathrm{N} A_j \Big|^{\frac{1}{p_0}}$$ for fixed $\mathrm{M}$ and $p_0 \ge 2$. In particular, it suffices to prove the reverse inequality. The idea is to show that the overlapping becomes small enough for large $\mathrm{M}$. Consider the partition $[-\Lambda,\Lambda]^n = \Phi_{\mathrm{M}, \Lambda} \cup \Psi_{\mathrm{M}, \Lambda}$ with $$\Psi_{\mathrm{M}, \Lambda} \ = \bigcup_{1 \le j \neq k \le \mathrm{N}} \Big( [-\Lambda, \Lambda]^n \cap \pi_\mathrm{M}^{g_j}(\gamma_{g_j^{-1}}(A_j)) \cap \pi_\mathrm{M}^{g_k}(\gamma_{g_k^{-1}}(A_k)) \Big),$$ the set where overlapping of the $\mathrm{B}_{\mathrm{M},j}$'s occur in $[-\Lambda,\Lambda]^n$. This yields
\begin{eqnarray*}
\lefteqn{\lim_{\mathrm{M} \to \infty} \mathrm{M}^n \Big\bracevert \Big( \sum_{j=1}^\mathrm{N} |\pi_\mathrm{M}^{g_j} \gamma_{g_j^{-1}} \chi_{A_j}|^2 \Big)^\frac12 \Big\bracevert_{p_0}^{p_0}} \\ & = & \lim_{\mathrm{M} \to \infty} \lim_{\Lambda \to \infty} \frac{\mathrm{M}^n}{2^n \Lambda^n} \int_{\Phi_{\mathrm{M},\Lambda}} \sum_{j=1}^\mathrm{N} \pi_\mathrm{M}^{g_j} \gamma_{g_j^{-1}} \chi_{A_j}(x) \, dx \\ & + & \lim_{\mathrm{M} \to \infty} \lim_{\Lambda \to \infty} \frac{\mathrm{M}^n}{2^n \Lambda^n}  \int_{\Psi_{\mathrm{M},\Lambda}} \Big( \sum_{j=1}^\mathrm{N} \pi_\mathrm{M}^{g_j} \gamma_{g_j^{-1}} \chi_{A_j}(x) \Big)^{\frac{p_0}{2}} \, dx \\ & \le & \sum_{j=1}^\mathrm{N} |A_j| \, + \, \mathrm{N}^{\frac{p_0}{2} - \frac12} \lim_{\mathrm{M},\Lambda \to \infty} \frac{\mathrm{M}^n}{2^n \Lambda^n} \int_{\Psi_{\mathrm{M},\Lambda}} \Big( \sum_{j=1}^\mathrm{N} \pi_\mathrm{M}^{g_j} \gamma_{g_j^{-1}} \chi_{A_j}(x) \Big)^{\frac{1}{2}} \, dx \\ & \le & \sum_{j=1}^\mathrm{N} |A_j| \, + \, \mathrm{N}^{\frac{p_0}{2} - \frac12} \lim_{\mathrm{M},\Lambda \to \infty} \frac{\mathrm{M}^n}{2^n \Lambda^n} \big| \Psi_{\mathrm{M},\Lambda} \big|^\frac12 \Big( \int_{[-\Lambda,\Lambda]^n} \sum_{j=1}^\mathrm{N} \pi_\mathrm{M}^{g_j} \gamma_{g_j^{-1}} \chi_{A_j}(x) \, dx \Big)^{\frac{1}{2}},
\end{eqnarray*}
where we have used the absence of overlapping in $\Phi_{\mathrm{M},\Lambda}$ to eliminate the power $p_0/2$ and $\mathrm{M}$-periodicity of the $\pi_\mathrm{M}^{g_j} \gamma_{g_j}^{-1} \chi_{A_j}$'s (isolatedly) for the first term. We may also use periodicity to estimate the last term on the right-hand side
\begin{eqnarray*}
\lefteqn{\lim_{\mathrm{M} \to \infty} \mathrm{M}^n \Big\bracevert \Big( \sum_{j=1}^\mathrm{N} |\pi_\mathrm{M}^{g_j} \gamma_{g_j^{-1}} \chi_{A_j}|^2 \Big)^\frac12 \Big\bracevert_{p_0}^{p_0}} \\ & \le & \sum_{j=1}^\mathrm{N} |A_j| \ + \ \mathrm{N}^{\frac{p_0}{2} - \frac12} \Big( \sum_{j=1}^\mathrm{N} |A_j| \Big)^\frac12 \Big( \lim_{\mathrm{M}, \Lambda \to \infty} \frac{\mathrm{M}^n}{2^n \Lambda^n} |\Psi_{\mathrm{M},\Lambda}| \Big)^\frac12 \\ & \le & \Big[ 1 + \mathrm{N}^{\frac{p_0}{2} - \frac12} \Big( \lim_{\mathrm{M} \to \infty} \lim_{\Lambda \to \infty} \frac{\mathrm{M}^n}{2^n \Lambda^n} |\Psi_{\mathrm{M},\Lambda}| \Big)^\frac12 \Big] \, \Big| \bigcup_{j=1}^\mathrm{N} A_j \Big|,
\end{eqnarray*}
since $\sum_j |A_j| \ge 1$ and $A_j \cap A_k = \emptyset$ for $j \neq k$. We write $\pi_\mathrm{M}^{\alpha_j}$ and $\pi_{\mathrm{M},\perp}^j$ for the $\mathrm{M}$-periodization operators in $\Theta_1$ and its orthocomplement respectively associated to $\rho_{g_j}$ (with rotation angle $\alpha_j$) and $\rho_{g_j,\perp}$. According to the form of $A_j$ we get the inequality
$$|\Psi_{\mathrm{M},\Lambda}| \le \sum_{1 \le j \neq k \le \mathrm{N}} |\Psi^1_{\mathrm{M},\Lambda}(j,k)| \times |\Psi^2_{\mathrm{M},\Lambda}(j,k)|$$ with the sets on the right given by
\begin{eqnarray*}
\Psi^1_{\mathrm{M},\Lambda}(j,k) & \!\! = \!\! & [-\Lambda,\Lambda]^2 \cap \pi_\mathrm{M}^{\alpha_j}(\rho_{g_j}^{-1}(R_j)) \cap \pi_\mathrm{M}^{\alpha_k}(\rho_{g_k}^{-1}(R_k)), \\ [3pt] \Psi^2_{\mathrm{M},\Lambda}(j,k) & \!\! = \!\! & [-\Lambda,\Lambda]^{n-2} \cap \pi_{\mathrm{M},\perp}^j \big( \rho_{g_j,\perp}^{-1}([-\lambda,\lambda]^{n-2}) \big) \cap \pi_{\mathrm{M},\perp}^k \big( \rho_{g_k,\perp}^{-1}([-\lambda,\lambda]^{n-2}) \big).
\end{eqnarray*}
Pick $z_0 \in \Theta_1$ and $\delta > 0$ minimal so that $\rho_{g_j}^{-1}(R_j) \subset \mathrm{B}_\delta(z_0)$ for all $j$.
By Lemma \ref{ErgodicLemma} and the hypothesis,
\begin{eqnarray*}
\lefteqn{ \hskip-10pt \lim_{\mathrm{M} \to \infty} \lim_{\Lambda \to \infty} \frac{\mathrm{M}^2}{4\Lambda^2} |\Psi^1_{\mathrm{M},\Lambda}(j,k)|} \\ & \le & \lim_{\mathrm{M} \to \infty} \lim_{\Lambda \to \infty} \frac{\mathrm{M}^2}{4\Lambda^2} \Big| [-\Lambda,\Lambda]^2 \cap \pi_\mathrm{M}^{\alpha_j}(\mathrm{B}_\delta(z_0)) \cap \pi_\mathrm{M}^{\alpha_k}(\mathrm{B}_\delta(z_0)) \Big| \ \lesssim \ \frac{1}{\mathrm{N}^{p_0+1} \lambda^{n-2}}.
\end{eqnarray*}
Moreover, eliminating the dependence on $k$ for $\Psi_{\mathrm{M},\Lambda}^2(j,k)$ gives the upper bound
\begin{eqnarray*}
\lefteqn{ \hskip-10pt \lim_{\mathrm{M} \to \infty} \lim_{\Lambda \to \infty} \frac{\mathrm{M}^{n-2}}{2^{n-2}\Lambda^{n-2}} |\Psi^2_{\mathrm{M},\Lambda}(j,k)|} \\ & \le & \lim_{\mathrm{M} \to \infty} \lim_{\Lambda \to \infty} \frac{\mathrm{M}^{n-2}}{2^{n-2}\Lambda^{n-2}} \Big| [-\Lambda,\Lambda]^{n-2} \cap \pi_{\mathrm{M},\perp}^j \big( \rho_{g_j,\perp}^{-1}([-\lambda,\lambda]^{n-2}) \big) \Big| \ \lesssim \ \lambda^{n-2}.
\end{eqnarray*}
Altogether, we have $\lim_{\mathrm{M}, \Lambda} \frac{\mathrm{M}^n}{2^n \Lambda^n} |\Psi_{\mathrm{M},\Lambda}| \lesssim \mathrm{N}^{1 - p_0}$ and the assertion follows.  \fin

\subsection{Infinite orbits admit Kakeya shadows} \label{PKaky}

We continue by introducing a class of sets in the unit sphere $\mathbb{S}^{n-1}$ for which the inequality in Lemma \ref{Decompactification} fails and show that infinite orbits of arbitrary discrete groups belong to this class. Our definition is motivated by Fefferman construction \cite{Fef}. Given a great circle $\zeta$ and any set of directions $\Omega$ in $\mathbb{S}^{n-1}$, we write $\pi_\zeta$ for the orthogonal projection onto the plane determined by $\zeta$. Let $$\Omega_{\zeta} \, = \, \Big\{ \frac{\pi_\zeta(\omega)}{|\pi_\zeta(\omega)|} \ \big| \ \omega \in \Omega \setminus \zeta^\perp \Big\} \, \subset \, \zeta$$ denote the geodesic projection of $\Omega$ onto $\zeta$. Let $\mathcal{R}_{\Omega_\zeta}$ denote the collection of rectangles $R$ in the plane determined by $\zeta$ with longest side pointing in a direction of $\Omega_{\zeta}$. The expression $3R$ will refer to the rectangle with the same center and width as $R$, but with $3$ times the length. We will say that $\Omega$ \emph{admits Kakeya shadows} if there exists an absolute constant $c_0$ such that for each $m \ge 1$ we may find a great circle $\zeta(m)$, a measurable set $E_m$ in the plane determined by $\zeta(m)$, and a finite collection of pairwise disjoint rectangles $\Sigma_{\Omega_{\zeta(m)}} \subset \mathcal{R}_{\Omega_{\zeta(m)}}$ so that
\begin{itemize}
\item[a)] $\displaystyle |E_m| \le \frac{c_0}{m} \sum_{R \in \Sigma_{\Omega_{\zeta(m)}}} |R|$,

\vskip3pt

\item[b)] $|R| \le c_0 |(3R \setminus R) \cap E_m|$ for each $R \in \Sigma_{\Omega_{\zeta(m)}}$.
\end{itemize}
The existence of sets admitting Kakeya shadows is a consequence of Besicovitch construction \cite{Be}. For instance, a set $\Omega$ admits Kakeya shadows whenever there exists a shadow $\Omega_\zeta$ which is dense in the unit circle $\mathbb{S}^1$. Let us now go back to the framework of our problem. Given a discrete group $\G$, a unit vector $u \in \R^n$ and an orthogonal representation $\gamma: \G \to O(n)$, we are interested in knowing when the $\gamma$-orbit of $u$ admits Kakeya shadows.

\begin{lemma} \label{Lema1}
The orbit $\mathcal{O}_\gamma(u)$ in $\mathbb{S}^{n-1}$ is either finite or admits Kakeya shadows.
\end{lemma}

\dem As described above, it suffices to prove the stronger statement that infinite orbits always admit a dense shadow.  The argument is simple when $n=2$. Assume the orbit $\mathcal{O}_\gamma(u)$ is not finite. By compactness, it must accumulate at some point $\sigma \in \mathbb{S}^1$ and we may find for each $\varepsilon > 0$ group elements $g_\varepsilon, h_\varepsilon \in \G$ such that $|\gamma_{g_\varepsilon}(u) - \gamma_{h_\varepsilon}(u)| < \varepsilon$ and $\mathrm{det} \gamma_{g_\varepsilon} = \mathrm{det} \gamma_{h_\varepsilon}$. It follows that $\gamma_{g_\varepsilon^{-1}h_\varepsilon}$ is a rotation of angle $< \varepsilon$. Density of $\mathcal{O}_\gamma(u)$ follows iterating these maps for $\varepsilon$ arbitrarily small. A similar argument applies when $\G$ is abelian and $n$ is arbitrary. In that case, the $\gamma_g$'s are pairwise commuting maps and we may consider again the direct sum decomposition into common eigenspaces $\R^n = \Delta_{\pm 1} \oplus \Theta_1 \oplus \Theta_2 \oplus \cdots \oplus \Theta_\ell$, where the $\gamma_g$'s act by rotations on two-dimensional $\Theta_j$'s. Let us decompose $u$ as $v_\delta + \sum_j v_j$ with $(v_\delta, v_j) \in \Delta_{\pm 1} \times \Theta_j$. If $\mathcal{O}_\gamma(u)$ is not finite, there must exist $1 \le j_0 \le \ell$ with $v_{j_0} \neq 0$ and $|\mathcal{O}_\gamma(v_{j_0})| = \infty$. Taking $\zeta_0$ to be the great circle in $\mathbb{S}^{n-1}$ generating $\Theta_{j_0}$, it is clear that the shadow $\mathcal{O}_\gamma(u)_{\zeta_0}$ contains $\mathcal{O}_\gamma(v_{j_0})$. Moreover, $\mathcal{O}_\gamma(v_{j_0})$ is an infinite orbit generated by rotations in $\Theta_{j_0}$, so that it is dense in $\zeta_0$ as in the $n=2$ case considered above.

Now that we know the statement holds for discrete abelian groups, it suffices to show that infinite orbits of discrete groups always admit infinite suborbits generated by an abelian subgroup. To see this, consider the linear subspace $$\mathcal{J}_\gamma \, = \, \Big\{ \xi \in \R^n \, \big| \, \gamma_{g^k}(\xi) = \xi \ \mbox{for all $g \in \G$ and some $k = k(g,\xi) \ge 1$} \Big\}.$$ If $u \notin \mathcal{J}_\gamma$, then we may find $g_0 \in \G$ such that $\gamma_{g_0}^j(u) \neq \gamma_{g_0}^k(u)$ for any pair of integers $j \neq k$. In that case, $\mathrm{A} = \langle g_0^k: \, k \in \Z \rangle \simeq \Z$ is an abelian subgroup of $\G$ giving rise to an infinite suborbit $\mathcal{O}_{\gamma}(\mathrm{A},u)$, as desired. On the other hand, if $u \in \mathcal{J}_\gamma$ we note that $\mathcal{J}_\gamma$ is a $\G$-module. Indeed, given a pair $(g_0,\xi_0) \in \G \times \mathcal{J}_\gamma$ and any $g \in \G$, define $h_0 = g_0^{-1} g g_0$ and $k_0 = k(h_0,\xi_0)$. Then we have by definition $$\gamma_{h_0^{k_0}}(\xi_0) = \xi_0 \, \Rightarrow \, \gamma_{g^{k_0}}(\gamma_{g_0}(\xi_0)) = \gamma_{g_0}(\xi_0).$$ This shows that $\gamma_{g_0}(\xi_0) \in \mathcal{J}_\gamma$ and proves our assertion above. In particular, since $u \in \mathcal{J}_\gamma$ the whole orbit $\mathcal{O}_\gamma(u)$ lies in $\mathcal{J}_\gamma$, which in turn is an invariant subspace of the representation $\gamma$. Restricting to that subspace if necessary, we may assume that $\mathcal{J}_\gamma = \R^n$. Let $e_1, e_2, \ldots, e_n$ denote an orthonormal basis of $\R^n$. Writing $m = m(g)$ for the least common multiple of $k(g,e_1), k(g,e_2), \ldots, k(g,e_n)$ it is easily checked that $\gamma_{g^m}(\xi) = \xi$ for all $\xi \in \R^n$. Let us consider the normal subgroup $\mathrm{H} = \{g \in \G \ | \ \gamma_g(\xi) = \xi \ \mbox{for all} \ \xi \}$. Again, since $\gamma_{\mid_\mathrm{H}}$ acts trivially on $\R^n$, we may restrict to the quotient group $\G/\mathrm{H}$ or equivalently assume that $\mathrm{H} = \{e\}$. This means that $g^m = e$ for all $g \in \G$ and some exponent $m = m(g) \ge 1$, so that $\G$ is a torsion group. According to the Jordan-Schur theorem, any torsion subgroup of the group of $n \times n$ complex matrices is virtually abelian. In other words, $\G$ admits an abelian normal subgroup $\mathrm{A}$ of finite order. If we now consider the suborbit $\mathcal{O}_{\gamma}(\mathrm{A},u)$, it must be infinite since $\mathcal{O}_\gamma(u)$ is infinite and $\mathrm{A}$ is of finite order. \fin

\begin{lemma} \label{Lema2}
Let $2 < p < \infty$, then $\mathcal{O}_\gamma(u)$ is finite if and only if \begin{eqnarray*}
\lefteqn{\hskip-10pt \Big\| \Big( \sum_{j=1}^\infty |H_{\gamma_{g_j}(u)}(f_{g_j})|^2 \Big)^\frac12 \Big\|_p} \\ & \le & c_p \max \left\{ \Big\| \Big( \sum_{j=1}^\infty |f_{g_j}|^2 \Big)^\frac12 \Big\|_p, \lim_{\mathrm{M} \to \infty} \mathrm{M}^{\frac{n}{p}} \Big\bracevert \Big( \sum_{j=1}^\infty |\pi_\mathrm{M}^{g_j} \gamma_{g_j^{-1}} f_{g_j}|^2 \Big)^\frac12 \Big\bracevert_p \right\},
\end{eqnarray*}
for any sequence $g_1, g_2, \ldots$ in $\G$ and any family of functions $f_{g_1}, f_{g_2}, \ldots$ in $L_p(\R^n)$.
\end{lemma}

\dem The validity of such an inequality for finite orbits and $1 < p < \infty$ follows from the $L_p$-boundedness of directional Hilbert transforms. Assume now that $\mathcal{O}_\gamma(u)$ is not finite and the inequality in the statement holds for some $p_0 > 2$. According to the proof of Lemma \ref{Lema1}, there exists an abelian subgroup $\mathrm{A}$ such that $\mathcal{O}_{\gamma}(\mathrm{A},u)$ is infinite. In other words, we may assume that $\G$ itself is abelian. Decomposing $\gamma$ into its irreducible parts as above, we get $\R^n = \Delta_{\pm 1} \oplus \Theta_1 \oplus \Theta_2 \oplus \cdots \oplus \Theta_\ell$ with $u = v_\delta + \sum_j v_j$. Pick $1 \le j_0 \le \ell$ such that $v_{j_0} \neq 0$ and $|\mathcal{O}_\gamma(v_{j_0})| = \infty$, and consider the great circle $\zeta_{j_0} = \Theta_{j_0} \cap \mathbb{S}^{n-1}$. We claim that there exists
\begin{itemize}
\item A family $\big\{ g_{km} \, | \, m \ge 1, \ 1 \le k \le 2^{2^m} \big\}$ in $\G$,

\vskip1pt

\item Rectangles $R_{km}$ and measurable sets $E_m$ in $\Theta_{j_0}$,
\end{itemize}
satisfying the following conditions
\begin{itemize}
\item[a)] $R_{km}$ points in the direction of $\pi_{\zeta_{j_0}}(\gamma_{g_{km}}(u))$,

\vskip3pt

\item[b)] $|E_m| \lesssim \frac{1}{m} \sum_k |R_{km}| \lesssim \frac{1}{m} \sum_k |(3R_{km} \setminus R_{km}) \cap E_m|$,

\vskip3pt

\item[c)] $R_{jm} \cap R_{km} = \emptyset = \gamma_{g_{jm}^{-1}}(R_{jm}) \cap \gamma_{g_{km}^{-1}}(R_{km})$ for $1 \le j \neq k \le 2^{2^m}$,

\item[d)] The sets $\{ (g_{km}, R_{km}) \}_{k}$ satisfy Lemma \ref{Ergodic2} with $\mathrm{N} = 2^{2^m}$ for each $m \ge 1$.
\end{itemize}
We finish the argument before proving the claim. Fix $m \ge 1$ in what follows and consider an orthonormal basis $e_1, e_2, \ldots, e_n$ so that $\mathrm{span} \{e_1,e_2\} = \Theta_{j_0}$. Let us write $\gamma_{g_{km}}(u)_j$ for $\langle \gamma_{g_{km}}(u), e_j \rangle$ and set $$\lambda \, = \, 1 + \max \left\{
\mathrm{length}(R_{km}) \frac{|\gamma_{g_{km}}(u)_j|}{\sqrt{\gamma_{g_{km}}(u)_1^2 + \gamma_{g_{km}}(u)_2^2}} \ \Big| \ 1 \le k \le 2^{2^m}, \ 3 \le j \le n \right\}.$$ We consider the following prisms in $\R^n$
\begin{eqnarray*}
A_{km} & = & \Big( R_{km} \times [-5\lambda,5\lambda]^{n-2} \Big), \\ B_{km} & = & \Big( (3R_{km} \setminus R_{km}) \times [-\mbox{$\frac12$} \lambda, \mbox{$\frac12$} \lambda]^{n-2} \Big).
\end{eqnarray*}
Taking $f_{g_{km}} = \chi_{A_{km}}$, we have $4|H_{\gamma_{g_{km}}(u)}(f_{g_{km}})| \ge \chi_{B_{km}}$. Indeed, by translational and rotational invariance we may assume that the prism $A_{km}$ is centered at $0$ and $\gamma_{g_{km}}(u)_1  = 0$. Now, given $x \in B_{km}$, we have
\begin{eqnarray*}
\lefteqn{\hskip-10pt \chi_{A_{km}}(x - t \gamma_{g_{km}}(u))} \\ & = & \chi_{[-\alpha, \alpha]}(x_1) \chi_{[-\beta, \beta]}(x_2 - t \gamma_{g_{km}}(u)_2) \prod_{j=3}^n \chi_{[-5\lambda, 5\lambda]}(x_j - t \gamma_{g_{km}}(u)_j)
\end{eqnarray*}
for some $0 < \alpha \le \beta = \mathrm{length}(R_{km})$. It is clear that $|x_1| \le \alpha$ for $x \in B_{km}$. On the other hand, the right hand side vanishes unless $|x_2 - t \gamma_{g_{km}}(u)_2| \le \beta$. This implies in turn that $$|t \gamma_{g_{km}}(u)_2| \le \beta + |x_2| \le 4\beta$$ $$\Downarrow$$ $$|x_j - t \gamma_{g_{km}}(u)_j| \le |x_j| + |t \gamma_{g_{km}}(u)_j| \le \frac12 \lambda + \frac{1}{\beta} |t \gamma_{g_{km}}(u)_2| \lambda \le \frac{9}{2} \lambda.$$ Hence, $\chi_{A_{km}}(x - t \gamma_{g_{km}}(u)) = \chi_{[-\beta,\beta]}(x_2 - t \gamma_{g_{km}}(u)_2)$ for $x \in B_{km}$ and we obtain $$H_{\gamma_{g_{km}}(u)}(f_{g_{km}})(x) = H(\chi_{[-\beta,\beta]})(x_2) = \frac{1}{\pi} \log \Big| \frac{x_2 + \beta}{x_2 - \beta} \Big| \quad \mbox{for} \quad x \in B_{km}.$$ The last logarithm is greater that $\frac{\pi}{4}$ for $\beta \le |x_2| \le 3\beta$, which gives the desired estimate $4|H_{\gamma_{g_{km}}(u)}(f_{g_{km}})| \ge \chi_{B_{km}}$. Using this, property b) above, and H\"older's inequality we obtain
\begin{eqnarray*}
\lambda^{n-2} \sum_{k=1}^{2^{2^m}} |R_{km}| & \lesssim & \lambda^{n-2} \sum_{k=1}^{2^{2^m}} |(3R_{km} \setminus R_{km}) \cap E_m| \\ & = & \sum_{k=1}^{2^{2^m}} \big| B_{km} \cap \big( E_m \times [-\mbox{$\frac12$}\lambda, \mbox{$\frac12$} \lambda]^{n-2} \big) \big| \\ & \lesssim & \int_{E_m} \int_{[-\frac12 \lambda, \frac12 \lambda]^{n-2}} \sum_{k=1}^{2^{2^m}} |H_{\gamma_{g_{km}}(u)}(f_{g_{km}})(x)|^2 \, dx \\ [5pt] & \le & |E_m|^{1 - \frac{2}{p_0}} \lambda^{(n-2)(1-\frac{2}{p_0})} \Big\|\Big( \sum_{k=1}^{2^{2^m}} |H_{\gamma_{g_{km}}(u)} (f_{g_{km}})|^2\Big)^\frac12 \Big\|^2_{p_0} \\ & \lesssim & \frac{1}{m^{1-\frac{2}{p_0}}} \Big( \lambda^{n-2} \sum_{k=1}^{2^{2^m}} |R_{km}| \Big)^{1 - \frac{2}{p_0}} \Big\|\Big( \sum_{k=1}^{2^{2^m}} |H_{\gamma_{g_{km}}(u)}(f_{g_{km}})|^2\Big)^\frac12 \Big\|^2_{p_0}.
\end{eqnarray*}
On the other hand, properties c) and d) give
\begin{eqnarray*}
\lefteqn{\hskip-10pt \max \left\{ \Big\| \Big( \sum_{k=1}^{2^{2^m}} |f_{g_{km}}|^2 \Big)^\frac12 \Big\|_{p_0}^2, \lim_{\mathrm{M} \to \infty} \mathrm{M}^{\frac{n}{p_0}} \Big\bracevert \Big( \sum_{k=1}^{2^{2^m}} |\pi_\mathrm{M}^{g_{km}} \gamma_{g_{km}^{-1}} f_{g_{km}}|^2 \Big)^\frac12 \Big\bracevert_{p_0}^2 \right\}} \\ & \lesssim & \Big( \sum_{k=1}^{2^{2^m}} \big| R_{km} \times [-5\lambda,5\lambda]^{n-2} \big| \Big)^{\frac{2}{p_0}} \ = \ (10)^\frac{2(n-2)}{p_0} \, \Big( \lambda^{n-2} \sum_{k=1}^{2^{2^m}} |R_{km}| \Big)^{\frac{2}{p_0}}
\end{eqnarray*}
with constants independent of $m$. Combining the two estimates via the inequality in the statement of the lemma, we get a contradiction for $m$ large enough. Therefore it suffices to prove our claim at the beginning of the proof. If $s = s(m) = 2^m$, let us write $F_m$ and $S_{1m}, S_{2m}, \ldots, S_{2^sm}$ for the measurable set $E$ and the pairwise disjoint rectangles $R_1, R_2, \ldots, R_{2^s}$ which appear in \cite[Lemma 10.1.1]{Gr}. Then we find $$|F_m| \lesssim \frac{1}{m} \summ_k |S_{km}| \lesssim \frac{1}{m} \summ_k |(3S_{km} \setminus S_{km}) \cap F_m|.$$ Of course, these properties remain stable under affine transformations and we may assume that $S_{1m}$ points in the direction $\pi_{\zeta_{j_0}}(u)$. Now, recall that the $\zeta_{j_0}$-shadow of $\mathcal{O}_\gamma(u)$ contains $\mathcal{O}_\gamma(v_{j_0})$, which in turn is dense in $\zeta_{j_0}$. Therefore, we may pick the group elements $g_{km}$ so that $\pi_{\zeta_{j_0}}(\gamma_{g_{km}}(u))$ points in the direction of $S_{km}$, or at least as close to it as we want. Note that Fefferman's construction admits small perturbations, so that we can move the rectangles slightly to make them point in the directions we pick. Our final choice of directions will be determined at the end of the proof. If we choose $R'_{km} \sim S_{km}$, it turns out that properties b) and c) hold except possibly for the pairwise disjointness condition $$\gamma_{g_{jm}^{-1}}(R'_{jm}) \cap \gamma_{g_{km}^{-1}}(R'_{km}) \, = \, \emptyset \quad \mbox{for} \quad 1 \le j \neq k \le 2^{2^m}.$$ Note however that the given construction is still invariant under translations. In particular, we may take $R'_{km} = \tau(S_{km})$ for a suitable translation $\tau$ in $\Theta_{j_0}$. Assume for clarity that $u = e_1$ and take $S_{1m} \sim [0,3 \log(s+2)] \times [0,2^{-s}]$, so that we roughly have $S_{km} \sim \gamma_{g_{km}}(S_{1m})$. This is essentially the worst possible configuration, since we find a large overlapping for the family $\gamma_{g_{km}^{-1}}(S_{km})$. Taking however $R'_{km} = \tau(S_{km})$ with $\tau(x,y) = (x+L,y)$ for $L>0$ large enough, it is easily checked that we get the missing disjointness condition in c). In summary, the rectangles $R_{km}'$ satisfy conditions b) and c) which in turn are stable under small perturbations. It remains to select the $g_{km}$'s and small perturbations $R_{km}$ so that conditions a) and d) also hold. Since the $g_{km}$'s are pairwise commuting because $\G$ can be taken abelian, it suffices to check that
\begin{itemize}
\item $R_{km} \parallel \pi_{\zeta_{j_0}}(\gamma_{g_{km}}(u))$,

\vskip12pt

\item $\displaystyle \summ_k |A_{km}| \ge 1$ for all $m \ge 1$,

\vskip4pt

\item $\big\langle \tan \alpha_{jk}^m \big\rangle \displaystyle \ge 2^{s(p_0+1)} \lambda^{n-2} \mathrm{diam}^2 \Big( \bigcup_{k=1}^{2^s} \gamma_{g_{km}}^{-1}(R_{km}) \Big) = J(p_0,m,\lambda)$.
\end{itemize}
Since $\lambda > 1$ by construction and $|R_{km}| \sim 2^{-s} \log s$, we see that $$\sum_{k=1}^{2^s} |A_{km}| \ge \sum_{k=1}^{2^s} |R_{km}| \gtrsim \log s \sim m >> 1.$$ Finally, we select the $g_{km}$'s so that the first and third conditions above hold. Pick $g_{1m} = e$, so that $R_{1m} = R_{1m}'$ points in the direction of $\pi_{\zeta_{j_0}}(u)$. Then we have to pick $g_{2m}$ so that $\langle \tan \alpha_{12}^m \rangle \ge J(p_0,m,\lambda)$ and $\pi_{\zeta_{j_0}}(\gamma_{g_{2m}}(u))$ is close enough to the direction of $S_{2m}$, so that $R_{2m}$ will be a small perturbation of $\tau(S_{2m})$ pointing in the direction of $\pi_{\zeta_{j_0}}(\gamma_{2m}(u))$. This is possible since the $\zeta_{j_0}$-shadow of $\mathcal{O}_\gamma(u)$ is dense in $\zeta_{j_0}$ and $\langle \tan \alpha_{12}^m \rangle \le J(p_0,m,\lambda)$ holds for finitely many directions. Once $g_{1m}$ and $g_{2m}$ are fixed, pick $g_{3m}$ with $\pi_{\zeta_{j_0}}(\gamma_{g_{3m}}(u))$ close enough to the direction of $S_{3m}$ and such that $\langle \tan \alpha_{j3}^m \rangle \ge J(p_0,m,\lambda)$ for $j=1,2$. Since the latter inequality fails for finitely many directions, again this is possible by density. Iterating the process we obtain the desired construction and the proof is complete. \fin

\demAA Since twisted Hilbert transforms are Fourier multipliers on the group von Neumann algebra associated to $\Gamma_{\mathrm{disc}}$, it is easily checked that $H_u \rtimes_\gamma id_\G$ is self-adjoint (up to conjugation of the symbol) so that we may assume $2 < p < \infty$. Then we combine Lemmas \ref{Decompactification} and \ref{Lema2}. \fin

\subsection{Twisted de Leeuw's compactification} \label{PdLw}

Let us now prove the implication i) $\Rightarrow$ ii) in Theorem A. We will write $H_u$ and $H_{\mathrm{b},u}$ in this paragraph to distinguish between the $u$-directional Hilbert transform on $\R^n$ and its Bohr compactification respectively. Assume $H_u \rtimes_\gamma id_\G$ is $L_p$-bounded for some $1 < p \neq 2 < \infty$. By density of trigonometric polynomials, it suffices to prove the $L_p$-boundedness of $H_{\mathrm{b},u} \rtimes_\gamma id_\G$ for finite sums of the form $$f_\mathrm{b} = \summ_{\xi,g} \widehat{f}(\xi,g) \hskip1pt \mbox{b-exp}_\xi \rtimes_\gamma \lambda_\G(g).$$ Once we have fixed $f_{\mathrm{b}}$, we find $(H_{\mathrm{b},u} \rtimes_\gamma id_\G)f_\mathrm{b} = (\widetilde{H}_{\mathrm{b},u} \rtimes_\gamma id_\G) f_\mathrm{b}$, where $$\widetilde{H}_{\mathrm{b},u}(\mbox{b-exp}_\xi) = -i \, \widetilde{\mathrm{sgn}} \langle u, \xi \rangle \mbox{b-exp}_\xi$$ and $\widetilde{\mathrm{sgn}}$ is a smoothing of the sign function which coincides with it on the finitely many values $\langle u, \xi \rangle$ with $\xi$ appearing in $f_\mathrm{b}$. We will write $\widetilde{H}_u$ for the corresponding smoothing operator in $\R^n$. Given $\delta > 0$, set $h_\delta(x) = (2 \pi \delta)^{-n/2} \exp(-|x|^2/2\delta)$ and we claim that
\begin{eqnarray*}
\lefteqn{\big\| (H_{\mathrm{b},u} \rtimes_\gamma id_\G) f_\mathrm{b} \big\|_{L_p(\widehat{\mathbf{\Gamma}_{\mathrm{disc}}})} = \big\| (\widetilde{H}_{\mathrm{b},u} \rtimes_\gamma id_\G) f_\mathrm{b} \big\|_{L_p(\widehat{\mathbf{\Gamma}_{\mathrm{disc}}})}} \\ [4pt] & = & \lim_{\delta \to \infty}  \Big\| -i h_\delta^{\frac{1}{p}} \summ_{\xi,g} \widetilde{\mathrm{sgn}} \langle u, \xi \rangle \widehat{f}(\xi,g) \exp_\xi \rtimes_\gamma \lambda_\G(g) \Big\|_{L_p(\widehat{\mathbf{\Gamma}})}  \\ [2pt] & = & \lim_{\delta \to \infty} \Big\| (\widetilde{H}_u \rtimes_\gamma id_\G) \Big(h_\delta^{\frac{1}{p}} \summ_{\xi,g} \widehat{f}(\xi,g) \exp_\xi \rtimes_\gamma \lambda_\G(g) \Big) \Big\|_{L_p(\widehat{\mathbf{\Gamma}})} \\ [2pt] & = & \lim_{\delta \to \infty} \Big\| (H_u \rtimes_\gamma id_\G) \Big(h_\delta^{\frac{1}{p}} \summ_{\xi,g} \widehat{f}(\xi,g) \exp_\xi \rtimes_\gamma \lambda_\G(g) \Big) \Big\|_{L_p(\widehat{\mathbf{\Gamma}})} \\ [2pt] & \lesssim & \lim_{\delta \to \infty} \Big\| h_\delta^{\frac{1}{p}} \summ_{\xi,g} \widehat{f}(\xi,g) \exp_\xi \rtimes_\gamma \lambda_\G(g) \Big\|_{L_p(\widehat{\mathbf{\Gamma}})} \ = \ \|f_{\mathrm{b}}\|_{L_p(\widehat{\mathbf{\Gamma}_{\mathrm{disc}}})}.
\end{eqnarray*}
Since the inequality follows by hypothesis, the constants are independent of the smoothing. It remains to justify the identities. If $f = \sum_{\xi,g} \widehat{f}(\xi,g) \exp_\xi \rtimes_\gamma \lambda_\G(g)$ stands for the trigonometric polynomial in $\mathcal{L}(\Gamma)$ with the same Fourier coefficients as $f_\mathrm{b}$, the second and fifth identities follow from $$\| f_\mathrm{b} \|_{L_p(\widehat{\mathbf{\Gamma}_{\mathrm{disc}}})} = \lim_{\delta \to \infty} \big\| h_\delta^{\frac{1}{p}} f \big\|_{L_p(\widehat{\mathbf{\Gamma}})}.$$ The proof reduces to $p=1$ since $h_\delta\rtimes_\gamma \mathbf{1}$ is in the center of $\mathcal{L}(\Gamma_\mathrm{disc})$, so $$\|f_\mathrm{b}\|_{L_p(\widehat{\mathbf{\Gamma}_{\mathrm{disc}}})}^p = \big\| |f_\mathrm{b}|^p \big\|_{L_1(\widehat{\mathbf{\Gamma}_{\mathrm{disc}}})} = \lim_{\delta \to \infty} \big\| h_\delta |f|^p \big\|_{L_1(\widehat{\mathbf{\Gamma}})} = \lim_{\delta \to \infty} \big\| h_\delta^{\frac{1}{p}} f \big\|_{L_p(\widehat{\mathbf{\Gamma}})}^p.$$ Approximating $|f_\mathrm{b}|^p$ by a finite sum, 
it suffices to note that $$\int_{\R^n_{\mathrm{bohr}}} \mbox{b-exp}_\xi(x) \, d\mu = \delta_{\xi=0} = \lim_{\delta \to \infty} \int_{\R^n} h_\delta(x) \exp_\xi(x) \, dx,$$ since the cross product with $\lambda_\G(g)$ only changes both sides by a factor of $\delta_{g=e}$. On the other hand, the third identity in our claim follows from the triangle inequality together with the identity $$\lim_{\delta \to \infty} \Big\| (\widetilde{H}_u \rtimes_\gamma id_\G) (h_\delta^{\frac{1}{p}} f) + i h_\delta^{\frac{1}{p}} \summ_{\xi,g} \widetilde{\mathrm{sgn}} \langle u, \xi \rangle \widehat{f}(\xi,g) \exp_\xi \rtimes_\gamma \lambda_\G(g) \Big\|_{L_p(\widehat{\mathbf{\Gamma}})} \, = \, 0.$$ Since the sums are finite, we prove that this is true term by term, in which case the cross products with $\lambda_\G(g)$ are irrelevant. In other words, we just need to show that we have $$\lim_{\delta \to \infty} \big\| \widetilde{H}_u(h_\delta^{\frac{1}{p}} \exp_\xi) + i h_\delta^{\frac{1}{p}} \widetilde{\mathrm{sgn}} \langle u, \xi \rangle \exp_\xi \big\|_p \, = \, 0.$$ According to the $L_p$-boundedness of $H_u$ in $\R^n$, these expressions are uniformly bounded in $\delta$ for fixed $1 < p < \infty$. By the three lines lemma, it then suffices to prove this identity in $L_2$ with exponents $(1-it)/2 + it/q = 1/2 + i\alpha$ with $\alpha \in \R$. By Plancherel theorem
\begin{eqnarray*}
\lefteqn{\big\| \widetilde{H}_u(h_\delta^{\frac{1}{2} + i \alpha} \exp_\xi) + i h_\delta^{\frac{1}{2}+i\alpha} \widetilde{\mathrm{sgn}} \langle u, \xi \rangle \exp_\xi \big\|_2^2} \\ & = & \int_{\R^n} \big| \widetilde{\mathrm{sgn}} \langle u, \eta \rangle - \widetilde{\mathrm{sgn}} \langle u, \xi \rangle \big|^2 \big| \widehat{h_\delta^{\frac12+i\alpha}}(\eta - \xi) \big|^2 \, d\eta \\ & \le &  \int_{|\xi - \eta| < \varepsilon} \big| \widetilde{\mathrm{sgn}} \langle u, \eta \rangle - \widetilde{\mathrm{sgn}} \langle u, \xi \rangle \big|^2 \big| \widehat{h_\delta^{\frac12+i\alpha}}(\eta - \xi) \big|^2 \, d\eta \\ & + & \int_{|\xi - \eta| \ge  \varepsilon} \big| \widetilde{\mathrm{sgn}} \langle u, \eta \rangle - \widetilde{\mathrm{sgn}} \langle u, \xi \rangle \big|^2 \big| \widehat{h_\delta^{\frac12+i\alpha}}(\eta - \xi) \big|^2 \, d\eta \ = \ A_{\delta,\varepsilon} + B_{\delta,\varepsilon}.
\end{eqnarray*}
Since $\widehat{h_\delta^{\frac12+i\alpha}}(\xi) = \delta^{\frac{n}{2}} \phi(\sqrt{\delta} \xi)$ for some Schwartz function $\phi$, we see that $B_{\delta,\varepsilon} \to 0$ as $\delta \to \infty$ for all $\varepsilon > 0$. On the other hand, $\widetilde{\mathrm{sgn}}$ is uniformly continuous and the integral of $h_\delta$ is $1$, so that $\sup_{\delta > 0} A_{\delta, \varepsilon} \to 0$ as $\varepsilon \to 0$. Combining both estimates we obtain the desired $0$ limit as $\delta \to \infty$. To justify the fourth identity we argue as for the third, so we may reduce it to show that $$\lim_{\delta \to \infty} \int_{\R^n} \big| \widetilde{\mathrm{sgn}} \langle u, \eta \rangle - \mathrm{sgn} \langle u, \eta \rangle \big|^2 \widehat{h_\delta^{{\frac12}+i \alpha}}(\eta - \xi)^2 \, d\eta \, = \, 0.$$ According to our smoothing, there exists $\varepsilon_0 > 0$ so that $\widetilde{\mathrm{sgn}} \langle u, \eta \rangle = \mathrm{sgn} \langle u, \eta \rangle$ if $|\eta - \xi| < \varepsilon_0$. In particular, the integral above is just defined on $|\eta - \xi| \ge \varepsilon_0$ and the limit again vanishes since $\mathrm{sgn}$ and $\widetilde{\mathrm{sgn}}$ are uniformly bounded functions. \fin

\subsection{Endpoint estimates for finite orbits} \label{Endpoint}

The proof of iii) $\Rightarrow$ i) in Theorem A is straightforward. We will instead prove stronger endpoint estimates from which we  recover i) by interpolation and duality, see \cite{JM2} for further details on $L_p$--$\mathrm{BMO}$ interpolation in this context. It is known that the directional Hilbert transform $H_u$ does not have a smooth Calder\'on-Zygmund kernel and fails to be $L_\infty \to \mathrm{BMO}$ bounded on $\R^n$ for $n > 1$ with the usual definition of $\mathrm{BMO}$. Nevertheless, there exists a \emph{directional $\mathrm{BMO}_u$ space} on $\R^n$ satisfying
\begin{itemize}
\item[a)] $H_u: L_\infty(\R^n) \to \mathrm{BMO}_u$,

\vskip1pt

\item[b)] $\big[ \mathrm{BMO}_u, L_p(\R^n) \big]_{p/q} = L_q(\R^n)$.
\end{itemize}
Namely, the norm in $\mathrm{BMO}_u$ is given by $$\|f\|_{\mathrm{BMO}_u} = \sup_{t > 0} \Big\| \Big( S_{u,t}|f|^2 - |S_{u,t}f|^2\Big)^{\frac12} \Big\|_\infty \quad \mbox{with} \quad \widehat{S_{u,t}f}(\xi) = e^{- t |\langle u,\xi \rangle |^2} \widehat{f}(\xi).$$ Properties a) and b) arose naturally in \cite{JMP1} by combining techniques from geometric group theory and diffusion semigroups, although it might be they were known previously to experts in the field. A similar result holds for $H_u \rtimes_\gamma id_\G$ when $|\mathcal{O}_\gamma(u)| < \infty$. Namely, setting $\widetilde{S}_{u,t} = S_{u,t} \rtimes_\gamma id_\G$ so that $$\widetilde{S}_{u,t}(\exp_\xi \rtimes_\gamma \lambda_\G(g)) \, = \, e^{-t |\langle u,\xi \rangle|^2} \exp_\xi \rtimes_\gamma \lambda_\G(g),$$ we define the space $\mathrm{BMO}_u(\widehat{\mathbf{\Gamma}})$ by $$\|f\|_{\mathrm{BMO}_u(\widehat{\mathbf{\Gamma}})} = \max \Big\{ \|f\|_{\mathrm{BMO}_u^r(\widehat{\mathbf{\Gamma}})}, \|f\|_{\mathrm{BMO}_u^c(\widehat{\mathbf{\Gamma}})} \Big\},$$ where the row and column BMO norms are given by
\begin{eqnarray*}
\|f\|_{\mathrm{BMO}_u^r(\widehat{\mathbf{\Gamma}})} & = & \Big\| \Big( \widetilde{S}_{u,t} (ff^*) - \widetilde{S}_{u,t}f \widetilde{S}_{u,t}f^* \Big)^\frac12\Big\|_{\mathcal{L}(\Gamma)}, \\ \|f\|_{\mathrm{BMO}_u^c(\widehat{\mathbf{\Gamma}})} & = & \Big\| \Big( \widetilde{S}_{u,t} (f^*f) - \widetilde{S}_{u,t}f^* \widetilde{S}_{u,t}f \Big)^\frac12\Big\|_{\mathcal{L}(\Gamma)}.
\end{eqnarray*}
See \cite{JM2} and the references therein for more on noncommutative BMO spaces.

\begin{theorem} \label{BMO}
If $\mathcal{O}_\gamma(u)$ is finite, then $$H_u \rtimes_\gamma id_\G: L_\infty(\widehat{\mathbf{\Gamma}}) \stackrel{cb}{\longrightarrow} \mathrm{BMO}_u(\widehat{\mathbf{\Gamma}}).$$ The superscript cb means that the map is not only bounded, but completely bounded.
\end{theorem}

\dem According to the $L_\infty \to \mathrm{BMO}$ analog of de Leeuw's compactification theorem ---already justified in \cite{JMP1}--- it suffices to show the complete boundedness of the map $$H_u \rtimes_\gamma id_\G: L_\infty (\widehat{\mathbf{\Gamma}_{\mathrm{disc}}}) \to \mathrm{BMO}_u(\widehat{\mathbf{\Gamma}_{\mathrm{disc}}}),$$ where the latter space is defined via the semigroup $$\mbox{b-exp}_\xi \rtimes_\gamma \lambda_\G(g) \mapsto e^{-t |\langle u,\xi \rangle|^2} \mbox{b-exp}_\xi \rtimes_\gamma \lambda_\G(g).$$ Let $\G_u = \big\{ g \in \G \, \big| \ \gamma_g(u) = u \big\}$ be the $\gamma$-stabilizer of $u$. Since the index $|\G : \G_u|$ coincides with $|\O_\gamma(u)|$, we have finitely many right cosets $\G_u g$. Let us label them as $\G_u g_j$ with $1 \le j \le |\O_\gamma(u)|$. Thus, if we set $$\Gamma_{\mathrm{disc}}^u = \R_{\mathrm{disc}}^n \rtimes_\gamma \G_u,$$ we may write any element $f$ in $\mathcal{L}(\Gamma_{\mathrm{disc}})$ as
\begin{eqnarray*}
f & = & \sum_{(\xi, g) \in \Gamma_{\mathrm{disc}}} \widehat{f}(\xi,g) \mbox{b-exp}_\xi \rtimes_\gamma \lambda_\G(g) \\ & = & \sum_{(\xi,g) \in \Gamma_{\mathrm{disc}}^u} \sum_{j=1}^{|\O_\gamma(u)|} \widehat{f}(\xi,gg_j) \mbox{b-exp}_\xi \rtimes_\gamma \lambda_\G(gg_j) \ = \ \sum_{j=1}^{|\O_\gamma(u)|} f_j.
\end{eqnarray*}
We define $$F_j = \sum_{(\xi,g) \in \Gamma_{\mathrm{disc}}} \widehat{f}(\xi,gg_j) \mbox{b-exp}_\xi \rtimes_\gamma \lambda_\G(g),$$ $w_j = \mbox{b-exp}_0 \rtimes_\gamma \lambda_\G(g_j)$, and let us also write $\mathsf{E}_u$ for the conditional expectation $\mathcal{L}(\Gamma_{\mathrm{disc}}) \to \mathcal{L}(\Gamma_{\mathrm{disc}}^u)$. Then, it is clear that we have $f = F_j w_j$ and $f_j = \mathsf{E}_u(F_j) w_j$ for all $j$, so that $\big( H_{u} \rtimes_\gamma id_{\G} \big) (f_j) = \big( H_{u} \rtimes_\gamma id_{\G_u} \big) (\mathsf{E}_u(F_j)) \, w_j$. Moreover, we have $\|\mathsf{E}_u(F_j)\|_\infty \le \|F_j\|_\infty = \|f\|_\infty$ since $w_j$ is a unitary. This yields
\begin{eqnarray*}
\big\| (H_{u} \rtimes_\gamma id_\G) (f) \big\|_{\mathrm{BMO}_u(\widehat{\mathbf{\Gamma}_{\mathrm{disc}}})} & \le & \sum_{j=1}^{|\O_\gamma(u)|} \big\| (H_{u} \rtimes_\gamma id_{\G}) (f_j) \big\|_{\mathrm{BMO}_u(\widehat{\mathbf{\Gamma}_{\mathrm{disc}}})} \\ & \le & |\O_\gamma(u)| \ \big\| H_{u} \rtimes_\gamma id_{\G_u} \big\|_{L_\infty \to \mathrm{BMO}} \ \|f\|_\infty.
\end{eqnarray*}
Since the same argument holds for matrix amplifications, it suffices to prove the cb-boundedness of $H_{u} \rtimes_\gamma id_{\G_u}$. In other words, we may assume that $u$ is a fixed point of $\gamma$. In that case we find
\begin{eqnarray*}
\widehat{H_{u}(\gamma_g f)}(\xi) & = & -i \hskip1pt \mathrm{sgn} \langle u , \xi \rangle \widehat{f}(\gamma_{g^{-1}}\xi) \\ & = & - i \hskip1pt \mathrm{sgn} \langle \gamma_g(u), \xi \rangle \widehat{f}(\gamma_{g^{-1}}\xi) \ = \ \widehat{\gamma_g(H_{u}f)}(\xi).
\end{eqnarray*}
This means that $H_{u}$ is a $\gamma$-equivariant map. On the other hand, since we know that $H_u: L_\infty(\R^n) \to \mathrm{BMO}_u(\R^n)$ the same holds for $H_{u}$ on $\R^n_{\mathrm{bohr}}$ ---use once more the $L_\infty \to \mathrm{BMO}$ analogue of de Leeuw's theorem from \cite{JMP1}--- and the assertion follows from a suitable application of the little Grothendieck inequality, see \cite{JMP1}. \fin

\begin{remark}
\emph{It follows from the proof that Theorem \ref{BMO} also holds for $\Gamma_\mathrm{disc}$.}
\end{remark}

\begin{theorem} \label{Weak}
If $\mathcal{O}_\gamma(u)$ is finite and $\G$ amenable, $$H_{u} \rtimes_\gamma id_\G: L_1(\widehat{\mathbf{\Gamma}_{\mathrm{disc}}}) \stackrel{cb}{\longrightarrow} L_{1,\infty}(\widehat{\mathbf{\Gamma}_{\mathrm{disc}}}).$$ 
\end{theorem}

\dem We may assume that $f$ is a positive trigonometric polynomial. Arguing as in the proof of Theorem \ref{BMO}, it suffices to consider the case where $u$ is fixed by $\gamma$. This means in particular that $(H_{u} \rtimes_\gamma id_\G)f$ is self-adjoint for $f$ positive. Let us consider the $*$-homomorphism $$\begin{array}{rrcl} \rho: & \mathcal{L}(\R^n_{\mathrm{disc}} \rtimes_\gamma \G) & \to & L_\infty(\R^n_{\mathrm{bohr}}) \bar\otimes \mathcal{B}(\ell_2(\G)) \\ & \mbox{b-exp}_\xi \rtimes_\gamma \lambda_\G(g) & \mapsto & \displaystyle  \summ_h \mbox{b-exp}_{\gamma_{gh}^{-1}(\xi)} \otimes e_{gh,h}. \end{array}$$ Since $\G$ is amenable, we may construct a F\o lner averaging sequence:
\begin{itemize}
\item $\G = \bigcup_{j \in J} \G_j$,

\vskip1pt

\item $|\G_j| < \infty$ and $\G_{j_1} \subset \G_{j_2}$ for $j_1 \le j_2$,

\vskip2pt

\item $|g\G_j \setminus \G_j| = o(|\G_j|)$ for all $g \in \G$.
\end{itemize}
If $p_j = \sum_{g \in \G_j} e_{g,g} \in \mathcal{B}(\ell_2(\G))$, we claim that following identities hold
\begin{eqnarray*}
\lefteqn{\big\| (H_{u} \rtimes_\gamma id_\G) f \big\|_{L_{1,\infty}(\widehat{\mathbf{\Gamma}_{\mathrm{disc}}})}} \\ & = & \sup_{\lambda > 0} \, \lambda \Big\|  \chi_{(\lambda, \infty)} \Big( \big| H_{u} \rtimes_\gamma id_\G) f \big| \Big) \Big\|_{L_1(\widehat{\mathbf{\Gamma}_{\mathrm{disc}}})} \\ & = & \sup_{\lambda > 0} \lambda \, \lim_j \frac{1}{|\G_j|} \Big\| p_j \rho \Big( \chi_{(\lambda, \infty)} \Big( \big| H_{u} \rtimes_\gamma id_\G) f \big| \Big) \Big) p_j \Big\|_{L_1(\R^n_{\mathrm{bohr}}; S_1(\ell_2(\G)))} \\ & = & \sup_{\lambda > 0} \lambda \, \lim_j \frac{1}{|\G_j|} \Big\| p_j \Big( \chi_{(\lambda, \infty)} \Big( \big| \rho \big[ ( H_{u} \rtimes_\gamma id_\G) f \big] \big| \Big) \Big) p_j \Big\|_{L_1(\R^n_{\mathrm{bohr}}; S_1(\ell_2(\G)))} \\ & = & \sup_{\lambda > 0} \lambda \, \lim_j \frac{1}{|\G_j|} \Big\| \chi_{(\lambda, \infty)} \Big( \big| p_j \rho \big[ ( H_{u} \rtimes_\gamma id_\G) f \big] p_j \big| \Big) \Big\|_{L_1(\R^n_{\mathrm{bohr}}; S_1(\ell_2(\G)))}.
\end{eqnarray*}
Let us first finish the argument assuming this claim. If $f = \sum_g f_g \rtimes_\gamma \lambda_\G(g)$,
\begin{eqnarray*}
\lefteqn{\hskip-10pt p_j \rho \big[ ( H_{u} \rtimes_\gamma id_\G) f \big] p_j} \\ & = & p_j \Big( \summ_{g,h} \gamma_{g^{-1}}(H_{u}(f_{gh^{-1}})) \otimes e_{g,h} \Big) p_j \\ & = & p_j \Big( \summ_{g,h} H_{u} (\gamma_{g^{-1}} (f_{gh^{-1}})) \otimes e_{g,h} \Big) p_j \\ & = & p_j \Big( \big[ H_{u} \otimes id_{\mathcal{B}(\ell_2(\G))} \big] (\rho f) \Big) p_j \ = \ \big[ H_{u} \otimes id_{\mathcal{B}(\ell_2(\G))} \big] (p_j \rho (f) p_j).
\end{eqnarray*}
The second identity uses that $u$ is fixed by $\gamma$. Writing Haar integration in $\R^n_{\mathrm{bohr}}$ as a limit of averages on arbitrary large cubes as we did in Paragraph \ref{Decompact}, a standard Fubini argument shows that the $L_1 \to L_{1,\infty}$ boundedness of $H_{u} \otimes id_{\mathcal{B}(\ell_2(\G))}$ reduces to that of $H \otimes id_{\mathcal{B}(\ell_2(\G))}$ where $H$ stands for the one-dimensional Hilbert transform in the Bohr compactification of $\R$. Such a weak type boundedness was proven by Randrianantoanina in his work \cite{R} on Hilbert transforms associated to maximal subdiagonal algebras. Thus, combining this with our claim we deduce that we have $$\big\| (H_{u} \rtimes_\gamma id_\G) f \big\|_{L_{1,\infty}(\widehat{\mathbf{\Gamma}_{\mathrm{disc}}})} \, \lesssim \, \lim_j \frac{1}{|\G_j|} \big\| p_j \rho (f) p_j \big\|_{L_1(\R^n_{\mathrm{bohr}}; S_1(\ell_2(\G)))} \, = \, \|f\|_{L_1(\widehat{\mathbf{\Gamma}_{\mathrm{disc}}})}.$$ The last identity above follows as in the second identity of our claim, which we now justify. The first identity is just the definition of the $L_{1,\infty}$ quasi-norm. The second follows from Neuwirth/Ricard's matrix-valued form of Szeg\"o's theorem \cite{NR}, some details ---also needed for the fourth identity--- can be found  below. The third follows since $\rho$ is a $*$-homomorphism. Indeed, $f$ is a trigonometric polynomial so that $(H_{u} \rtimes_\gamma id_\G) f$ is a bounded operator. Therefore, we may replace $\chi_{(\lambda,\infty)}$ by $\chi_{(\lambda, \mathrm{M})}$ for $\mathrm{M}$ large enough and argue by polynomial approximation. 
The last identity can be proved by following \cite[Proof of Theorem 2.1]{NR} again, where the idea is to approximate $\chi_{(\lambda,\mathrm{M})} | \cdot |$ by polynomials and estimate $$\lim_j \frac{1}{|\G_j|} \Big\| p_j P \Big( \rho \big[ (H_{u} \rtimes_\gamma id_\G)f \big] \Big) p_j - P \Big( p_j \rho \big[ (H_{u} \rtimes_\gamma id_\G)f \big] p_j \Big) \Big\|_1$$ for each polynomial $P$. As we have $$p_j x^k p_j - (p_j x p_j)^k = p_k x^{k-1} (x p_j - p_j x p_j) + (p_j x^{k-1} p_j - (p_j x p_j)^{k-1}) x p_j,$$ an induction argument yields the inequality $$\big\| p_j x^k p_j - (p_j x p_j)^k \big\|_1 \, \le \, (k-1) \|x\|_\infty^{k-1} \big\| x p_j - p_j x p_j \big\|_1.$$ On the other hand, analyzing the trigonometric polynomial $f$ term by term, we are reduced to showing that $|\G_j|^{-1} \|A p_j - p_j A p_j\|_1 \to 0$ for $A = \rho (\mbox{b-exp}_\xi \rtimes_\gamma \lambda_\G(g))$ and this follows from the relation $$\big\| A p_j - p_j A p_j \big\|_1 = \Big\| \sum_{h \in \G_j \setminus g^{-1} \G_j} \mbox{b-exp}_{\gamma_{gh}^{-1}(\xi)} \otimes e_{gh,h} \Big\|_1 \le |g \G_j \setminus \G_j|$$ and the fact that  $|g\G_j \setminus \G_j| = o(|\G_j|)$, which follows from the amenability of $\G$. \fin

\demA We have proved i) $\Leftrightarrow$ ii) $\Leftrightarrow$ iii) and the endpoint estimates. The remaining equivalence with iv) is now very simple. Indeed, the boundedness for finite orbits follows from the triangle inequality since any block of rows is contractively complemented in the Schatten $p$-class $S_p(\G)$. When the orbit is infinite, unboundedness follows by picking $f_{g,h} = \delta_{h=e} f_g$ so that the resulting square function inequality only holds for orbits with no Kakeya shadows. \fin

\begin{remark}
\emph{Given a Fourier multiplier $T_m$ on $\R^n$ ---for example, the directional Hilbert transform in this paper---  and a orthogonal representation $\gamma: \G \to O(n)$, we may consider three noncommutative forms of such an operator}
\begin{itemize}
\item[(A)] \emph{The matrix operator $f_{g,h} \otimes e_{g,h} \mapsto \gamma_{g^{-1}}T_m(f_{g,h}) \otimes e_{g,h}$,}

\item[(B)] \emph{The cross product operator $f_g \rtimes_\gamma \lambda_\G(g) \mapsto T_m(f_g) \rtimes_\gamma \lambda_\G(g)$,}

\item[(C)] \emph{The cocycle form of the multiplier in $\mathcal{L}(\G)$: $\lambda_\G(g) \mapsto m_{b(g)} \lambda_\G(g)$.}
\end{itemize}
\emph{In terms of $L_p$-boundedness, \cite{NR} gives (A) $\Rightarrow$ (B) for $\G$ discrete amenable and  \cite{JMP1} ---see Paragraph \ref{Last} below--- gives (B) $\Rightarrow$ (C) for arbitrary discrete $\G$. One could wonder when the reverse implications hold. When dealing with directional Hilbert transforms, Theorem A shows (A) $\Leftrightarrow$ (B) for any discrete $\G$, while the comment after Corollary C proves that the implication (C) $\Rightarrow$ (B) fails in general.}
\end{remark}

\begin{remark}
\emph{Given the unboundedness for infinite orbits and our $L_\infty \to \mathrm{BMO}$ estimate for finite orbits ---which hold in the category of operator spaces--- we see that $L_p$-boundedness is equivalent to complete $L_p$-boundedness for twisted Hilbert transforms and $1 < p < \infty$, which was not clear \emph{a priori}.}
\end{remark}

\section{Lacunarity, cocycles, and convergence of Fourier series}

In this section we analyze the more general frameworks considered in Theorem B and Corollary C. We will also establish some connections between these problems and the $L_p$-boundedness of directional maximal operators or idempotent Fourier multipliers on $\R$.

\subsection{Lacunary $\gamma$-suborbits}

Given $n \ge 2$ and a set of directions $\Omega \subset \mathbb{S}^{n-1}$ in the unit sphere, the directional maximal operator $M_\Omega$ is defined on smooth functions $f: \R^n \to \C$ by $$M_\Omega f(x) = \sup_{w \in \Omega} \sup_{r > 0} \frac{1}{2r} \int_{-r}^r |f(x - t \omega)| \, dt.$$ The sets $\Omega$ in the circle for which $M_\Omega$ is bounded in $\R^2$ can now be described with remarkable accuracy. Bateman recently proved in \cite{B} that $M_\Omega$ is $L_q$-bounded for some/any $1 < q < \infty$ iff $\Omega$ is a finite union of lacunary sets of finite order in the sense of Sj\"ogren/Sj\"olin \cite{SS}. In higher dimensions, the only known results were due to Carbery and Nagel/Stein/Wainger \cite{C,NSW}. In a recent paper \cite{PR}, we obtain more general results and characterize the $L_q$-boundedness for arbitrary dimensions. In particular, we prove that $M_\Omega$ is $L_q$-bounded for any $1 < q < \infty$ provided that $\Omega$ is $\mathrm{HD}$-lacunary (see the introduction for the definition). This will be the key ingredient in the following lemma.

\begin{lemma} \label{SqMax}
If $1 < p < \infty$ and $\Omega$ is $\mathrm{HD}$-lacunary, then $$\Big\| \Big( \sum_{\omega \in \Omega} |H_\omega f_\omega|^2 \Big)^\frac12 \Big\|_p \, \lesssim \, \Big\| \Big( \sum_{\omega \in \Omega} |f_\omega|^2 \Big)^\frac12 \Big\|_p.$$
\end{lemma}

\dem Since $H_\omega$ is essentially self-dual and the case $p=2$ is clear, we may clearly assume that $p > 2$. Let $\frac1q = 1 - \frac2p$, then we select v in the positive part of the unit ball of $L_q(\R^n)$ such that
\begin{eqnarray*}
\Big\| \Big( \sum_{\omega \in \Omega} |H_\omega f_\omega|^2 \Big)^\frac12 \Big\|_p^2 & = & \Big\| \sum_{\omega \in \Omega} |H_{\omega} f_\omega|^2 \Big\|_{\frac{p}{2}} \\ & = & \sum_{\omega \in \Omega} \int_{\R^n} |H_{\omega}f_\omega|^2(x) \, \mathrm{v}(x) \, dx \\ & = & \sum_{\omega \in \Omega} \int_{\R^n \ominus \R \omega} \Big( \int_{\R \omega} \big| H(f_{z,\omega})(s) \big|^2 \mathrm{v}_{z,\omega}(s) \, ds \Big) dz,
\end{eqnarray*}
where $f_{z,\omega}(s) = f_\omega (z + s \omega)$ and $\mathrm{v}_{z,\omega}(s) = \mathrm{v}(z + s \omega)$ for $z \perp \omega$ and $H$ stands for the Hilbert transform on $\R$. Now we pick some $0 < \delta < 1$ and use the Hardy-Littlewood maximal operator $M$ on $\R$ to get $$\mathrm{v}_{z,\omega}(s) \le M(\mathrm{v}_{z,\omega}^{\frac{1}{\delta}})^\delta(s) = \sup_{r > 0} \Big( \frac{1}{2r} \int_{-r}^r \mathrm{v}^{\frac{1}{\delta}} \big( z+ (s-t) \omega \big) \, dt \Big)^\delta \le M_\Omega (\mathrm{v}^\frac{1}{\delta} )^\delta(z+s\omega).$$ It is well-known \cite{Duo} that $\mathrm{w}_\delta = M(\mathrm{v}_{z,\omega}^{1/\delta})^\delta$ is an $A_2$ Muckenhoupt weight with $A_2$ contants depending only on $\delta$. Since the Hilbert transform is bounded on $L_2(\R, \mathrm{w}_\delta(s) ds)$ with norm depending ---linearly, see \cite{Pet}--- on the $A_2$ norm of $\mathrm{w}_\delta$ we conclude
$$\Big\| \sum_{\omega \in \Omega} |H_\omega f_\omega|^2 \Big\|_{\frac{p}{2}} \lesssim \sum_{\omega \in \Omega} \int_{\R^n} | f_{\omega}(x) |^2 M_\Omega(\mathrm{v}^{\frac{1}{\delta}})^\delta(x) \, dx \le \|M_\Omega\|_{q\delta \to q\delta}^\delta \, \Big\|\sum_{\omega \in \Omega} |f_\omega|^2 \Big\|_{\frac{p}{2}}.$$ This also follows from \cite{CF2}. We now use HD-lacunarity and the result from \cite{PR}. \fin

\begin{remark} \label{psi-lac}
\emph{Given a discrete group $\G$ and a length function $\psi: \G \to \R_+$, a countable subset $\Delta = \{\delta_j \, | \, j \ge 1\} \subset \G$ will be called \emph{$\psi$-lacunary} when the following condition holds $$\sup_{j \ge 1} \frac{\psi(\delta_{j+1})}{\psi(\delta_j)} \, \le \, \lambda_\psi \, < \, 1.$$ Arguing as in Lemma \ref{LPReduction1}, we may construct a sequence of smooth functions $h_m$ on $\R_+$ fulfilling the hypotheses of Lemma \ref{Lema3} for the group $\Gamma_\mathrm{disc}$ and the length function $\xi \rtimes_\gamma g \mapsto \psi(g)$, so that $$\|f\|_{L_p(\widehat{\mathbf{\Gamma}_{\mathrm{disc}}})} \sim \Big\| \sum_{\delta \in \Delta} f_\delta \rtimes_\gamma \lambda_\G(\delta) \otimes e_\delta \Big\|_{L_p(\widehat{\mathbf{\Gamma}_\mathrm{disc}}; \ell_{rc}^2)}$$ for every $f \in L_{\Delta,p}(\widehat{\mathbf{\Gamma}_\mathrm{disc}})$ ($1 < p < \infty$) and constants depending only on $p,\lambda_\psi$.}
\end{remark}

\demB Enumerating $\Lambda = \{g_j \, | \, j \ge 1\}$, we know by hypothesis that there exists some $\omega$ in the unit sphere so that $\gamma_{g_j}^{-1}(u) \to \omega$ lacunarly as $j \to \infty$. In particular, if $1 \le j \le \mathrm{M} << \mathrm{N}$ we see that
\begin{eqnarray*}
\big| \gamma_{g_{j+1}^{-1}}(u) - \gamma_{g_\mathrm{N}^{-1}}(u) \big| & \le & \big| \gamma_{g_\mathrm{N}^{-1}}(u) - \omega \big| + \big| \gamma_{g_{j+1}^{-1}}(u) - \omega \big| \\ & \le & (1+\lambda) \big| \gamma_{g_\mathrm{N}^{-1}}(u) - \omega \big| + \lambda \big| \gamma_{g_j^{-1}}(u) - \gamma_{g_\mathrm{N}^{-1}}(u) \big|.
\end{eqnarray*}
In particular, if $\mathrm{N} = \mathrm{N(M)}$ is large enough we find for $1 \le j \le \mathrm{M}$
\begin{eqnarray*}
\big| \gamma_{g_{j+1}g_\mathrm{N}^{-1}}(u) - u \big| & = & \big| \gamma_{g_{j+1}^{-1}}(u) - \gamma_{g_\mathrm{N}^{-1}}(u) \big| \\ & \le & \sqrt{\lambda} \, \big| \gamma_{g_{j}^{-1}}(u) - \gamma_{g_\mathrm{N}^{-1}}(u) \big| \ = \ \sqrt{\lambda} \, \big| \gamma_{g_j g_\mathrm{N}^{-1}}(u) - u \big|.
\end{eqnarray*}
For sufficiently large $\mathrm{M}$ we approximate the $p$-norm of $f \in L_{\Lambda,p}(\widehat{\mathbf{\Gamma}_\mathrm{disc}})$
$$\|f\|_{L_p(\widehat{\mathbf{\Gamma}_{\mathrm{disc}}})} \, \sim \, \Big\| \sum_{j=1}^\mathrm{M} f_{g_j} \rtimes_\gamma \lambda_\G(g_j) \Big\|_{L_p(\widehat{\mathbf{\Gamma}_{\mathrm{disc}}})} \, = \, \Big\| \sum_{j=1}^\mathrm{M} f_{g_j} \rtimes_\gamma \lambda_\G(g_j g_\mathrm{N}^{-1}) \Big\|_{L_p(\widehat{\mathbf{\Gamma}_{\mathrm{disc}}})}.$$
We may now apply Remark \ref{psi-lac}. Indeed, consider
\begin{itemize}
\item The set $\Delta_{\mathrm{M}} = \{ g_j g_\mathrm{N}^{-1} \, | \, 1 \le j \le \mathrm{M} \}$,

\item The length function $\psi_{\gamma,u}(g) = |\gamma_g(u) - u|^2$.
\end{itemize}
We refer to Paragraph \ref{PLength-gamma} to justify that $\psi_{\gamma,u}$ is a length. According to our estimates above, we see that $\Delta_{\mathrm{M}}$ is $\psi_{\gamma,u}$-lacunary and Remark \ref{psi-lac} yields the following norm equivalence with $\delta_j = g_j g_\mathrm{N}^{-1}$ $$\|f\|_{L_p(\widehat{\mathbf{\Gamma}_{\mathrm{disc}}})} \sim \Big\| \sum_{\delta \in \Delta_{\mathrm{M}}} f_{\delta g_\mathrm{N}} \rtimes_\gamma \lambda_\G(\delta) \otimes e_\delta \Big\|_{L_p(\widehat{\mathbf{\Gamma}_{\mathrm{disc}}}; \ell_{rc}^2)}.$$ Using the same equivalence for $(H_u \rtimes_\gamma id_\G)f$, we are reduced to proving $$\Big\| \sum_{\delta \in \Delta_{\mathrm{M}}} \!\! H_u(f_{\delta g_\mathrm{N}}) \rtimes_\gamma \lambda_\G(\delta) \otimes e_\delta \Big\|_{L_p(\widehat{\mathbf{\Gamma}_\mathrm{disc}}; \ell_{rc}^2)} \, \lesssim \, \Big\| \sum_{\delta \in \Delta_{\mathrm{M}}} \!\! f_{\delta g_\mathrm{N}} \rtimes_\gamma \lambda_\G(\delta) \otimes e_\delta \Big\|_{L_p(\widehat{\mathbf{\Gamma}_\mathrm{disc}}; \ell_{rc}^2)}$$
for $1 < p < \infty$ and constants independent of $\mathrm{M}$. As explained in Section \ref{Sect1}, these norms are sums/intersections of row and column spaces
for $p$ smaller/greater than $2$. In particular, it suffices to show that this inequality holds for row and column spaces independently.  In the row case, the inequality reads as $$\Big\| \Big( \sum_{j=1}^\mathrm{M} \big| H_u(f_{g_j}) \big|^2 \Big)^\frac12 \Big\|_p \, \lesssim \, \Big\| \Big( \sum_{j=1}^\mathrm{M} |f_{g_j}|^2 \Big)^\frac12 \Big\|_p$$ which clearly holds from the $L_p$-bdness of $H_u$. In the column case we have $$\Big\| \Big( \sum_{j=1}^\mathrm{M} \big| H_{\gamma_{g_\mathrm{N} g_j^{-1}}(u)}(\gamma_{g_\mathrm{N} g_j^{-1}} f_{g_j}) \big|^2 \Big)^\frac12 \Big\|_p  \, \lesssim \, \Big\| \Big( \sum_{j=1}^\mathrm{M} |\gamma_{g_\mathrm{N} g_j^{-1}}f_{g_j}|^2 \Big)^\frac12 \Big\|_p$$ for functions $f_g \in L_p(\R^n_{\mathrm{bohr}})$. Arguing as in Paragraph \ref{PdLw}, we are reduced to proving such an inequality in the Euclidean space $L_p(\R^n)$. Now we use our second assumption which gives $\mathrm{HD}$-lacunarity for the suborbit $\mathcal{O}_\gamma(\Lambda^{-1},u)$ ---and therefore also for the set $\gamma_{g_\mathrm{N}} \gamma_{g_j}^{-1}(u)$--- in conjunction with Lemma \ref{SqMax} to deduce the validity of such a square function inequality in $L_p(\R^n)$ for $1 < p < \infty$ with absolute constants independent of $\mathrm{M}$. This completes the proof. \fin

\begin{remark} \label{GralLac}
\emph{Theorem B admits several generalizations. Namely, we could work with other length functions $\psi$ for which $\Lambda$ were $\psi$-lacunary as long as $\Lambda$ is a finite covering of the suborbit $\mathcal{O}_\gamma(\Lambda,u)$: $\sup_{g \in \Lambda} \big| \big\{ h \in \Lambda \, | \, \gamma_h(u) = \gamma_g(u) \big\} \big| \, < \, \infty$. On the other hand, more general notions of $\psi$-lacunarity may be considered. It would be interesting to obtain Littlewood-Paley estimates for $\psi$-lacunary sequence of higher order in the sense of \cite{SS}, with which one could relax the conditions in Theorem B.}
\end{remark}

\begin{remark} \label{RMeyer}
\emph{We have found in Lemma \ref{LPReduction1} a twisted form on $L_p(\R^n_\mathrm{bohr})$ of Meyer's square function inequality. Now we may provide necessary and sufficient conditions for this inequality to hold. Indeed, it follows from the proof of Theorem A that not admitting Kakeya shadows is necessary for the orbits/suborbits considered. On the other hand, being $\mathrm{HD}$-lacunary is sufficient, as we see from Lemma \ref{SqMax} and de Leeuw's compactification like in Paragraph \ref{PdLw}. Simpler arguments ---ergodicity and transference $\R^n \leftrightarrow \R^n_{\mathrm{bohr}}$ are not needed in the Euclidean-Lebesguean case--- yield the same conclusions for the twisted Meyer's inequality on $L_p(\R^n)$. It is an open problem to decide if there exist sets of directions which do not admit Kakeya shadows and which fail to be $\mathrm{HD}$-lacunary.}
\end{remark}

\subsection{Convergence of Fourier series in the $\psi$-metric} \label{Last}

Corollary C follows from our results above on the discretized algebra $\mathcal{L}(\Gamma_\mathrm{disc})$, an intertwining identity from \cite{JMP1} and standard Fourier methods.

\demC If $\dim \H = n$ the mapping $\pi_\psi: \mathcal{L}(\G) \to \mathcal{L}(\R_{\mathrm{disc}}^n) \rtimes_{\gamma} \G$ determined by $\lambda_\G(g) \mapsto \mbox{b-exp}_{b(g)} \rtimes_{\gamma} \lambda_\G(g)$ is a trace preserving $*$-homomorphism. The key property is that $$\pi_\psi \circ H_{\psi,u} \, = \, (H_u \rtimes_{\gamma} id_\G) \circ \pi_\psi,$$ which can be easily checked. Note also that $\pi_\psi(L_{\Lambda,p}(\widehat{\mathbf{G}})) = L_{\Lambda,p}(\widehat{\mathbf{\Gamma}_\mathrm{disc}})$. This allows us to represent $H_{\psi,u}$ as the restriction of the $\gamma$-twisted Hilbert transform to the image of $\pi_\psi$. In particular, the assertions in a) and b) on the boundedness of $H_{\psi,u}$ follow from the corresponding boundedness of $H_u \rtimes_{\gamma} id_\G$ considered in Theorems A and B. On the other hand, the $L_p$-density of trigonometric polynomials ---for which the convergence results hold trivially--- allows us to emulate the standard argument in $\mathbb{T}^n$ for $L_p$ convergence of Fourier series. In other words, we must show that $$f \mapsto \sum_{g \in \G} \chi_{\mathrm{RK}}(b(g)) \widehat{f}(g) \lambda_\G(g)$$ defines an $L_p$-bounded Fourier multiplier with constants independent of $\mathrm{R}$. By the intertwining identity above, it suffices to prove uniform $L_p$-boundedness for the family $T_{\mathrm{RK}} \rtimes_{\gamma} id_\G$, where $T_{\mathrm{RK}}$ is the Fourier multiplier in $L_p(\R^n_{\mathrm{bohr}})$ with Fourier symbol $\chi_{\mathrm{RK}}$. If we denote the faces of $\mathrm{K}$ by $\partial_j \mathrm{K}$ ($1 \le j \le m$), this in turn factorizes as a finite product of semispace Fourier multipliers of the form $S_{u_j, v_j} \rtimes_{\gamma} id_\G$ with $u_j \perp \partial_j \mathrm{K}$, $v_j \in \mathrm{R} \partial_j \mathrm{K}$ and $$\widehat{S_{u_j, v_j}f}(\xi) = \chi_{\R_+} \langle \xi - v_j, u_j \rangle \widehat{f}(\xi) \quad \Rightarrow \quad S_{u_j,v_j}f = M_{v_j} \circ S_{u_j} \circ M_{-v_j} f$$ where $S_uf = \frac12 (id + i H_u)$ and $M_{v}f = \mbox{b-exp}_v f$. Since the modulations $M_v \rtimes_{\gamma} id_\G$ are $L_p$-isometries, the convergence result in c) follows once again from Theorem A applied to $H_{u_j} \rtimes_{\gamma} id_\G$. It remains to justify the necessity in c). Consider the cocycle $(\mathbf{H}, \mathbf{b}, \boldsymbol\gamma)$ in $\Gamma_{\mathrm{disc}}$ defined as follows $$\mathbf{H} = \R^n_{\mathrm{disc}}, \qquad \mathbf{b}(\xi, g) = \xi, \qquad \boldsymbol\gamma(\xi, g) = \gamma_g.$$ The associated $\boldsymbol\psi(\xi,g) = |\xi|^2$ yields $H_{\boldsymbol\psi,u} = H_u \rtimes_\gamma id_\G$ and finiteness of $\gamma(\G)$ gives $$\lim_{\mathrm{R} \to \infty} \Big\| f - \sum_{\xi \in \mathrm{RK}}^{\null} \sum_{g \in \G} \widehat{f}(\xi,g) \lambda_{\Gamma_{\mathrm{disc}}}(\xi \rtimes_\gamma g) \Big\|_{L_p(\widehat{\mathbf{\Gamma}_{\mathrm{disc}}})} = 0$$ from the same argument above. Now we assume that the limit above vanishes. By a standard application of the uniform boundedness principle,  we deduce that $\sup_{\mathrm{R} > 0} \| T_{\mathrm{RK}} \rtimes_\gamma id_\G \|_{p \to p} < \infty$. Since we have already seen how translations of the Fourier symbol can be written in terms of conjugation against isometric modulation maps, we also deduce that we must have $$\sup_{\mathrm{R} > 0} \big\| T_{\tau_{\mathrm{R},j} \mathrm{RK}} \rtimes_\gamma id_\G \big\|_{p \to p} \, < \, \infty,$$ where $T_{\tau_{\mathrm{R},j} \mathrm{RK}}$ is the Fourier multiplier associated to the symbol $\chi_{\tau_{\mathrm{R},j} \mathrm{RK}}$. If we pick the translations $\tau_{\mathrm{R},j}$ so that $\partial_j \mathrm{K} \cap \partial_j \tau_{\mathrm{R},j} \mathrm{RK} \neq \emptyset$, we see that $\tau_{\mathrm{R},j} \mathrm{RK}$ approximates the semispace determined by the face $\partial_j \mathrm{K}$. Applying Fatou's lemma, we conclude that $H_{u_j} \rtimes_\gamma id_\G$ must be $L_p$-bounded for all $1 \le j \le m$. The result finally follows from another application of Theorem A. \fin

\begin{remark}
\emph{As with Remark \ref{GralLac}, we may also consider different lengths for c).}
\end{remark}


\begin{remark}
\emph{As it was justified in the introduction, condition a) in Corollary C no longer provides a characterization of the $L_p$-boundedness of $H_{\psi,u}$ for arbitrary discrete groups. Such a characterization appears to be harder and would yield examples of idempotent Fourier multipliers on group von Neumann algebras for arbitrary discrete groups. A characterization of the $L_p$ boundedness of $H_{\psi,u}$ for $\G = \R$ ---i.e. idempotent Fourier multipliers determined by restriction from inner cocycles--- will appear in \cite{PPR}.}
\end{remark}

\noindent \textbf{Acknowledgement.} We thank Andrei Jaikin-Zapirain for discussions regarding the geometric form of infinite orbits on the sphere and especially Marius Junge for allowing us to reproduce the cross product generalization \cite{JMP2} of de Leeuw's compactification. 
Javier Parcet has been supported in part by ERC Starting Grant StG-256997-CZOSQP. Keith Rogers was supported in part by ERC Starting Grant StG-277778-RESTRICTION. Both authors were supported in part by the Spanish grant MTM2010-16518 and ICMAT Severo Ochoa project SEV-2011-0087.

\bibliographystyle{amsplain}


\vskip15pt

\hfill \noindent \textbf{Javier Parcet} \\ \null \hfill\texttt{javier.parcet@icmat.es} \\
\null \hfill Instituto de Ciencias Matem{\'a}ticas \\ \null \hfill
CSIC-UAM-UC3M-UCM \\ \null \hfill Consejo Superior de
Investigaciones Cient{\'\i}ficas \\ \null \hfill C/ Nicol\'as Cabrera 13-15.
28049, Madrid. Spain

\vskip10pt
	
\hfill \noindent \textbf{Keith Rogers} \\ \null \hfill\texttt{keith.rogers@icmat.es}
\\ \null \hfill Instituto de Ciencias Matem{\'a}ticas \\ \null \hfill
CSIC-UAM-UC3M-UCM \\ \null \hfill Consejo Superior de
Investigaciones Cient{\'\i}ficas \\ \null \hfill C/ Nicol\'as Cabrera 13-15.
28049, Madrid. Spain

\end{document}